\documentclass[12pt]{amsart}
\usepackage{amssymb,latexsym,amsmath,amscd}
\theoremstyle{plain}
\newtheorem{thm}{Theorem}[section]
\newtheorem{prop}[thm]{Proposition}
\newtheorem{cor}[thm]{Corollary}
\newtheorem{lemma}[thm]{Lemma}
\theoremstyle{definition}
\newtheorem{defn}[thm]{Definition}
\newtheorem{rem}[thm]{Remark}
\newtheorem{notat}[thm]{Notation}
\newtheorem{ex}[thm]{Example}

\newcommand{\ra}{\rightarrow}
\newcommand{\lra}{\longrightarrow}
\newcommand{\PP}{\mathbf{P}}

\newcommand{\ZZ}{\mathbf{Z}}

\newcommand{\Q}{\mathcal{Q}}
\newcommand{\bin}[2]{ {{#1} \choose {#2}} }
\newcommand{\OO}{\mathcal{O}}

\newcommand{\HH}{\mathbf{H}}
\newcommand{\I}{\mathcal{I}}
\newcommand{\M}{\mathcal{M}}
\newcommand{\MM}{\mathfrak{m}}

\newcommand{\FF}{\mathbf{F}}
\newcommand{\GG}{\mathbf{G}}

\newcommand{\OP}[1]{\OO_{\PP^{#1}}}

\newcommand{\mmax}{{\rm max}}
\newcommand{\mmin}{{\rm min}}

\title{The General Hyperplane Section of a Curve}

\author{Elisa Gorla}

\address{Department of Mathematics
\\ University of Notre Dame \hfil\break\indent 255 Hurley Hall, Notre Dame, IN
46556-4618, USA}

\email{egorla@nd.edu}

\thanks{The author was partially supported by a scholarship from the Italian research institute ``Istituto
Nazionale di Alta Matematica Francesco Severi''}

\begin{document}

\maketitle

Abstract: In this paper, we discuss some necessary and sufficient conditions for a curve
to be arithmetically Cohen-Macaulay, in terms of its general hyperplane section.
We obtain a characterization of the degree matrices that can occur for points in
the plane that are the general plane section of a non arithmetically Cohen-Macaulay curve of $\PP^3$.
We prove that almost all the degree matrices with positive subdiagonal that occur for the general plane
section of a non arithmetically Cohen-Macaulay curve of $\PP^3$, arise also as degree matrices of some
smooth, integral, non arithmetically Cohen-Macaulay curve, and we characterize the exceptions. 
We give a necessary condition on the graded Betti
numbers of the general plane section of an arithmetically Buchsbaum (non arithmetically Cohen-Macaulay)
curve in $\PP^n$. For curves in $\PP^3$, we show that any set of Betti numbers that satisfy that condition
can be realized as the Betti numbers of the general plane section of an arithmetically Buchsbaum, non
arithmetically Cohen-Macaulay curve. We also show that the matrices that arise as degree matrix of the
general plane section of an arithmetically Buchsbaum, integral, (smooth) non arithmetically
Cohen-Macaulay space curve are exactly those that arise as degree matrix of the general plane section of
an arithmetically Buchsbaum, non arithmetically Cohen-Macaulay space curve and have positive
subdiagonal. 
We also prove some bounds on the dimension of the deficiency module of an arithmetically Buchsbaum space
curve in terms of the degree matrix of the general plane section of the curve, and we prove that they are 
sharp.

\newpage

It is well known that several invariants of an arithmetically Cohen-Macaulay projective scheme, such as
the degree, the $h$-vector, the graded Betti numbers, and many more, are preserved when we intersect
the scheme with a hyperplane that meets it properly. Moreover, the intersection of an arithmetically
Cohen-Macaulay scheme of dimension at least $1$ with a hyperplane is itself arithmetically Cohen-Macaulay.
If we are interested in a $d$-dimensional, arithmetically Cohen-Macaulay scheme $V\subset\PP^n$, we can
intersect it with a hyperplane that meets it properly. Repeating the procedure $d$ times, we get a
zero-dimensional scheme $X\subset\PP^{n-d}$. Then we can deduce the invariants of $V$ from the invariants of $X$.

In the general case of a scheme that is not necessarily arithmetically Cohen-Macaulay, not even all the
hyperplane sections will have the same invariants. However, a generic hyperplane $H$ will intersect
$V$ properly, and the scheme $V\cap H$ will always have the same invariants. In general, though, the
invariants of $V$ cannot be easily deduced from those of the general hyperplane section $V\cap H$. 
In the case when $V\cap H$ is arithmetically Cohen-Macaulay and has dimension at least $1$, however, $V$
itself is forced to be arithmetically Cohen-Macaulay. In particular, we are again in the situation when
we can deduce invariants of $V$, from those of $V\cap H$.

A great deal of work has been devoted to the analysis of the case when $V\cap H$ has dimension $0$, or
equivalently when $V$ is
a projective curve. Obviously, we cannot expect to deduce the Cohen-Macaulayness of $V$ from the
Cohen-Macaulayness of $X$, with no further assumptions. In fact, the general hyperplane section of a curve
is a zero-dimensional scheme, so it is always arithmetically Cohen-Macaulay. A.V. Geramita and J.C.
Migliore, R. Strano, R. Re, C. Huneke and B. Ulrich, found sufficient conditions on the general
hyperplane section of a curve, that guarantee Cohen-Macaulayness of the curve (see \cite{GM2}, \cite{S},
\cite{R}, \cite{HU}, \cite{M2}). 
A brief summary and discussion of the work that has been done in the papers we just mentioned is contained 
in {\em Section 1} of this paper. Section 1 contains some terminology and
notation as well. We also introduce the concept of lifting matrix of a zero-dimensional scheme $X\subset\PP^n$ (see
Definition~\ref{liftingmatrix}). The lifting matrix is a matrix of integers, whose entries are the differences between
the shifts of the last and first free module in a minimal free resolution of $X$.

The starting point of {\em Section 2} is a sufficient condition found by C.~Huneke and B.~Ulrich
for $V$ to be arithmetically Cohen-Macaulay, in terms of the graded Betti numbers of its general
hyperplane section (see Theorem~\ref{socle}, Corollary~\ref{bigger3} and Corollary~\ref{bigger3P^n}). 
For example, for a curve in~$\PP^3$ the general plane section is a zero-dimensional scheme $X$ in~$\PP^2$. 
The matrix of integers whose entries are the degrees of the entries of the Hilbert-Burch matrix of~$X$
is called the degree matrix of $X$. A sufficient condition for the curve to be arithmetically
Cohen-Macaulay, is that all the entries of the degree matrix of~$X$ are at least~$3$. 
The question we want to answer is: is this condition necessary as well? That is, can we construct a non
arithmetically Cohen-Macaulay curve, whose general plane section has a prescribed degree matrix, for each
degree matrix that has at least one entry less than or equal to $2$? In Theorem~\ref{main} and
Theorem~\ref{mmain}, we prove that the sufficient condition of Huneke and Ulrich is necessary as well. We
do so by constructing a non
arithmetically Cohen-Macaulay curve, whose general plane section has a prescribed degree matrix, for each
degree matrix that has one entry less than or equal to~$2$. 
The curves we construct in Theorem~\ref{main} are connected and reduced, and they are the
union of two arithmetically Cohen-Macaulay curves. The construction of Theorem~\ref{main},
however, requires a further assumption on one of the entries of the degree matrix of $X$, in case it has
size bigger than $2\times 3$. The curves we construct in Theorem~\ref{mmain} are a union of smooth,
connected complete intersections. The construction of Theorem~\ref{mmain} works in full generality, for
any degree matrix that has one of the entries smaller than or equal to $2$. Moreover, we ask whether it is
possible to give a necessary condition for Cohen-Macaulayness of such a curve, in terms of the
$h$-vector of its general plane section. As one can expect, the answer to this question is
negative, as we show in Proposition~\ref{vect}.

In {\em Section 3} we deal with integral (that is, reduced and irreducible) curves in $\PP^3$. 
We ask whether it is possible to
find a condition on the degree matrix of the general plane section of a curve, which is weaker
than assuming that all the entries are bigger than or equal to~$3$, but still forces Cohen-Macaulayness of
the curve under the hypothesis that the curve is integral. Moreover, we ask whether it is possible to 
give a sufficient condition for Cohen-Macaulayness of an integral curve
in~$\PP^3$, in terms of the $h$-vector of its general plane section.
We are able to produce two families of degree matrices that do not have all the
entries bigger than or equal to~$3$, but with the property that any integral curve whose general plane 
section has one of those degree matrices is arithmetically Cohen-Macaulay.
So we have sufficient conditions on the degree matrix of the general plane section of a
curve that, together with integrality of the curve, force the curve to be arithmetically Cohen-Macaulay.
They are treated in Proposition~\ref{aCM2x3} and Proposition~\ref{aCM-II}.  
From those, we are able to deduce sufficient conditions for Cohen-Macaulayness of an integral curve, 
in terms of the $h$-vector of its general plane section. In particular, the curve has that same
$h$-vector. This is shown in Corollary~\ref{h-vect}. 
In Theorem~\ref{smooth2x3} and Theorem~\ref{main-smooth}, we show that, except for the two families
treated in Proposition~\ref{aCM2x3} and Proposition~\ref{aCM-II}, the degree matrices with positive
subdiagonal that correspond to points that are the general plane section of a non arithmetically
Cohen-Macaulay curve, are the same as the degree matrices that correspond to points that are the general
plane section of a non arithmetically Cohen-Macaulay, integral curve. Notice that the degree matrix of a
zero-dimensional scheme that is the general plane section of an integral curve needs to have positive entries 
on the subdiagonal. For each degree matrix that does not fall in the two categories of
Proposition~\ref{aCM2x3} and Proposition~\ref{aCM-II}, we construct a smooth, connected, non
arithmetically Cohen-Macaulay curve, whose general plane section has that degree matrix. It follows that
any admissible $h$-vector of decreasing type, except for those treated in Corollary~\ref{h-vect},
can be realized as the
$h$-vector of the general plane section of an integral, (or even smooth and connected) non arithmetically 
Cohen-Macaulay curve. This is proven in Corollary~\ref{hvects}. Notice that any admissible
$h$-vector of decreasing type can be realized as the $h$-vector of an integral, arithmetically
Cohen-Macaulay curve in~$\PP^3$, hence of its general plane section (this follows for example from 
\cite{HTV}).

In {\em Section 4} we concentrate on arithmetically Buchsbaum curves in $\PP^{n+1}$. We investigate 
whether we can give some conditions on the Betti numbers of the general plane section of an 
arithmetically Buchsbaum, non arithmetically Cohen-Macaulay curve.
In Proposition~\ref{BuchP^n}, we look at the lifting matrix (defined in Section 1)
of a zero-dimensional scheme which is the general hyperplane section of an arithmetically Buchsbaum, 
non arithmetically
Cohen-Macaulay curve. We show that one of the entries of such a lifting matrix has to be equal to $n$. 
For the case of curves in~$\PP^3$, the lifting matrix of the general plane section coincides with its
degree matrix. Therefore, the degree matrix of the general plane section of an arithmetically Buchsbaum, 
non arithmetically Cohen-Macaulay curve in $\PP^3$ has at least one entry equal to~$2$.
In Theorem~\ref{mainBuch} we show that this condition is both necessary and
sufficient. We do so by constructing an
arithmetically Buchsbaum curve whose general plane section has a prescribed degree matrix, for any
degree matrix that has at least one entry equal to~$2$. Then we analyze the case of integral,
arithmetically Buchsbaum, non arithmetically Cohen-Macaulay curves of~$\PP^3$. The general plane section
of an integral curve is a set of points in Uniform Position, hence its degree matrix has positive
subdiagonal. In Theorem~\ref{intBuch}, we show that for any degree matrix whose subdiagonal is positive,
and that has at least one entry equal to~$2$, 
we can construct a smooth, connected, arithmetically Buchsbaum, non arithmetically Cohen-Macaulay
curve in~$\PP^3$, whose general plane section has that degree matrix. In other words, a homogeneous
matrix of integers occurs as degree matrix of the general plane section of some integral,
(or smooth and connected) arithmetically Buchsbaum, non arithmetically Cohen-Macaulay curve if and only 
if it has positive subdiagonal and at least one entry is equal to~$2$.
We also prove some bounds for the dimension of the deficiency module of an arithmetically
Buchsbaum curve, degree by degree in Proposition~\ref{BuchBounds}, and globally in Corollary~\ref{dimMB}.
The bounds are again in terms of the entries of the degree matrix of the general plane section of the
curve. In the end of Section~4, we produce families of examples of arithmetically Buchsbaum, non
arithmetically Cohen-Macaulay curves of~$\PP^3$ that achieve the previously mentioned bounds, in order to
show their sharpness. The curves that we produce have general plane section that is either level or whose
homogeneous saturated ideal is generated in a single degree.

The author would like to thank J. Migliore for many useful discussions.
The computer algebra system CoCoA (\cite{CoCoA}) was used during the preparation of this paper, to compute
examples and verify some statements.

\section{Preliminaries and Notation}

Let $C$ be a curve in $\PP^{n+1}=\PP^{n+1}(k)$, where $k$ is an algebraically
closed field. In Section~3 and part of Section~2, we will assume that $k$ has characteristic~0. Throughout
the paper, a curve will be a non-degenerate,  equidimensional, locally  Cohen-Macaulay, dimension~1
subscheme of~$\PP^{n+1}$.  

Let $I_C$ be the saturated homogeneous ideal corresponding to $C$ in the 
polynomial ring $S=k[x_0,x_1,\ldots,x_{n+1}]$. We will denote by $\MM$ the
homogeneous irrelevant maximal ideal of $S$, $\MM=(x_0, x_1, \ldots ,x_{n+1})$.
Let $\I_C\subset \OP{n+1}$ be the ideal sheaf of~$C$.

We will denote the \emph{cohomology modules of $C$} by
$$H^i_*(\I_C)=\bigoplus_{m\in\ZZ} H^i(\PP^{n+1},\I_C(m))$$
and will denote the dimension of their graded pieces as $$h^i(\I_C(m))=dim_k\; H^i(\I_C(m)).$$ 
The first cohomology module of a curve $C$ is also called \emph{deficiency module}. We will denote it by 
$\M_C$.

\begin{notat}
For $M$ an $R$-module, we denote by $\alpha(M)$ the
initial degree of the module 
$$\alpha(M)=\mmin \{ m\in \ZZ\; |\; M_m\neq 0 \}.$$
If $M$ has finite length, we denote by $\alpha^+(M)$ its final degree
$$\alpha^+(M)=\mmax \{ m\in \ZZ\; |\; M_m\neq 0 \}.$$
\end{notat}

It is well known that the deficiency module of $C$ is trivial if and only if $C$ is
arithmetically Cohen-Macaulay. Its deficiency module has finite length as an $S$-module (or
equivalently, finite dimension as a $k$-vector space) if and only if $C$ is locally Cohen-Macaulay and
equidimensional (see \cite{S}, \cite{HIO} 37.5 or \cite{M}, Theorem 1.2.5).

In this paper, we will extend a result of R.~Strano (\cite{S}) and a result of C.~Huneke and 
B.~Ulrich (\cite{HU}). We are interested in finding conditions on the general hyperplane
section  of $C$, that are necessary and sufficient for the Cohen-Macaulayness of the curve.
$X$ will denote the zero-dimensional scheme that is the general hyperplane section of $C$ and 
$I_X$ its homogeneous, saturated ideal in the polynomial ring $R=k[x_0,x_1,\ldots,x_n]$. Sometimes, we
will also use $I_X$ for the ideal of $X$ as a subset of $\PP^{n+1}$, i.e. $I_X$ will be an ideal of
$S=k[x_0,x_1,\ldots,x_{n+1}]$.

We will devote particular attention to space curves $C\subset\PP^3$. In this case,
$I_X$ is a codimension 2 Cohen-Macaulay ideal of $R=k[x_0,x_1,x_2]$, hence a standard 
determinantal ideal, due to the Hilbert-Burch Theorem (see \cite{E}, Theorem 20.15). 
It is generated by the maximal minors of a $t\times (t+1)$ homogeneous matrix $A=(F_{ij})$. Let
$M=(a_{i,j})$ be the matrix whose entries are the degrees of the entries of $A$; $M$ is the {\it degree
matrix} of $X$. We make the convention that the entries of $M$ decrease from right to left and from top
to bottom: $a_{i,j} \leq a_{k,r}$, if $i\geq k$ and $j \leq r$. If some entry $F_{ij}$ 
of $A$ is $0$, then the degree is not well defined. In this case, there exist $k,l$ such that
$F_{ik}, F_{lk}, F_{lj}$ are all different from zero. We set $a_{ij}=a_{ik}-a_{lk}+a_{lj}.$

We can assume without loss of generality that the Hilbert-Burch matrix
has the property that $F_{ij}=0$ if $a_{ij}\le 0$, and
$deg(F_{ij})=a_{ij}$ if $a_{ij}>0$. Note that some of the $F_{ij}$'s could be $0$ even if $a_{ij}>0$.

A matrix of integers $M=(a_{i,j})$ is {\em homogeneous} if 
$a_{i,j}+a_{r,s}=a_{i,s}+a_{r,j}$ for all $i,r=1,\ldots ,t$ and $j,s=1,\ldots ,t+1$.
Notice that the degree matrix of a homogeneous matrix is homogeneous in this sense.
Abusing language, we use the term {\em degree matrix} to refer to any matrix of integers that 
is the degree matrix of some scheme in projective space.

A {\em standard determinantal} scheme $X\subseteq\PP^n$ of codimension $c$, is a scheme whose saturated
ideal $I_X$ is minimally generated by the maximal minors of a matrix of polynomials of size $t\times
(t+c-1)$, for some $t$. The definition of standard determinantal scheme was introduced by M.~Kreuzer,
J.C.~Migliore, U.~Nagel and C.~Peterson in \cite{KMNP}. In particular, any
Cohen-Macaulay ideal of codimension~2 is standard determinantal. We can characterize the matrices of
integers that are also degree matrices of some standard determinantal scheme, as those that are
homogeneous and whose diagonal is entirely positive.

\begin{prop}\label{diag_entries}
Let $M=(a_{i,j})$ be a matrix of integers of size $t\times (t+c-1)$. Then $M$ is a degree matrix if and
only if it is homogeneous and $a_{h,h}>0$ for $h=1,\ldots,t$.
\end{prop} 

\begin{proof}
Any degree matrix is homogeneous, as we observed before. 
We start by showing that every degree matrix has positive entries on the diagonal. 
We will prove the thesis by contradiction, showing that if $a_{h,h}\leq 0$ for some~$h$, 
then the scheme~$X$ cannot be standard determinantal. So let 
$A=(F_{i,j})_{i=1,\ldots,t; \; j=1,\ldots,t+c-1}$ be the
matrix defining~$X$; equivalently, $I_X$ is minimally generated by the maximal minors of $A$. 
In particular, the
determinant~$\Delta$ of the submatrix $B=(F_{i,j})_{i=1,\ldots,t; \; j=1,\ldots,t}$ is nonzero.
Assume $a_{h,h}\leq 0$ for some $h$, then $a_{i,j}\leq 0$ for $i\geq h$ and $j\leq h$. Hence $F_{i,j}=0$
for $i\geq h$ and $j\leq h$. Then~$B$ contains a submatrix of zeroes of size $(t-h+1)\times h$. 

We claim that $\Delta=0$. Let us prove it by induction on the size $t$ of $B$.
For $t=1$, we have $B=(0)$; $B$ is a matrix of size $1\times 1$.
Assume now that the thesis is true for $t-1$ and prove it for $t$. We have
$$\Delta=\sum_{i=1}^{t}(-1)^{i+t}F_{i,t}b_{i,t}$$
where $b_{i,j}=det(B_{i,j})$ is the determinant of the submatrix $B_{i,j}$, obtained from $B$ deleting
the $i$-th row and the $j$-th column. For each $i$, $B_{i,t}$ is a matrix of size $(t-1)\times(t-1)$ that
has a submatrix of $h$ columns and (at least) $t-h$ rows consisting of zeroes. Thus, induction hypothesis
applies on $B_{i,t}$ for all $i$, giving $b_{i,t}=0$. So $\Delta=0$, contradicting the
assumption that $X$ is standard determinantal.

Conversely, let $M=(a_{i,j})$ be a homogeneous matrix of integers of size $t\times (t+c-1)$, 
with positive diagonal. We want to show that $M$ is a
degree matrix. We need to exhibit a standard determinantal scheme that has $M$ as its degree matrix.
So let $$A=\left(\begin{array}{cccccc}
F_{1,1} & \cdots & F_{1,c} & 0 & 0 & \cdots \\
0 & F_{2,2} & \cdots & F_{2,c+1} & 0 & \cdots \\
 & & \ddots & & \ddots & \\
0 & 0 & \cdots & F_{t,t} & \cdots & F_{t,t+c-1} \end{array}\right)$$
where $F_{i,j}\in R$ are generic homogeneous polynomials of degree $deg(F_{i,j})=a_{i,j}$. By assumption,
all the degrees involved are positive.
$A$ defines a standard determinantal, reduced scheme (see \cite{BCG}, Proposition 2.5), whose saturated
homogeneous ideal is minimally generated by the maximal minors of $A$.
\end{proof}

Let us consider the general case of curves embedded in a projective space of arbitrary dimension.
If $C\subset \PP^{n+1}$, a general hyperplane section of $C$ is a zero-dimensional scheme $X\subset\PP^n$. 
We would like to associate a matrix of integers to each zero-dimensional scheme, such that it extends
the idea of degree matrix to arbitrary codimension.

\begin{defn}\label{liftingmatrix}
Let $X\subset\PP^n$ be a zero-dimensional scheme with minimal free resolution
$$0\lra\FF_n=\bigoplus_{i=1}^{t}R(-m_i)\lra\FF_{n-1}\lra\cdots\lra\FF_2\lra\FF_1=
\bigoplus_{j=1}^{r}R(-d_j)\lra I_X\lra 0$$
where $m_1\geq\ldots\geq m_t$ and $d_1\geq\ldots\geq d_r$. 

The matrix of integers $M=(a_{ij})=(m_i-d_j)$ is the \emph{lifting matrix} of $X$.
\end{defn} 

Notice that the lifting matrix coincides with the degree matrix of $X$ in the case of space curves ($n=2$). 
The lifting matrix will play the role of the degree matrix of $X$, for $n>2$. Notice moreover,
that the entries of $M$ decrease from right to left and from top
to bottom: $a_{i,j} \leq a_{k,l}$, if $i\geq k$ and $j\leq l$.

A {\em complete intersection} of {\em type} $(d_1,\ldots,d_r)$ is a scheme whose
homogeneous,  saturated ideal is generated by a regular sequence of forms of degrees $d_1\leq
d_2\leq\ldots\leq d_r$. We will abbreviate it by $CI(d_1,\ldots,d_r)$, or by $CI$ when we do not need to
specify the degrees.

We will always assume that the curve $C\subset\PP^{n+1}$ is non-degenerate. Notice that for $n=2$, if $C$
is degenerate then it is a plane curve, so it is arithmetically Cohen-Macaulay. We will often
abbreviate {\em arithmetically Cohen-Macaulay} by {\em aCM}.

We can assume that if the zero-dimensional scheme $X\subset\PP^n$is the general hyperplane section of a 
non-degenerate $C\subset\PP^{n+1}$, then $X$ is non-degenerate, as the following Lemma shows. 
See \cite{HU}, or Proposition 2.2 in \cite{M2} for a proof. The Lemma extends a result of O. A. Laudal
for curves in $\PP^3$ (see \cite{L}, pg. 142 and 147) . 

\begin{lemma}
The general hyperplane section of a non-degenerate curve $C\subset\PP^{n+1}$ of degree $d\geq n+1$ is
non-degenerate.
\end{lemma}

The case $t=1$, $n=2$, that is the case when the general plane section of $C\subset\PP^3$ is a complete 
intersection, has been studied by R. Strano. He proved the following result (Theorem 6, \cite{S}).

\begin{thm}\label{strano}
Let $C\subset\PP^3$ be a reduced and irreducible, non-degenerate curve of degree $d$ not lying on
a quadric surface. If the general plane section $X$ is a $CI(s,t)$, then $C$ is a $CI(s,t)$.
\end{thm}

The result is sharp, in the sense that we can easily find examples of curves that are non-aCM, whose 
general plane section is a complete intersection of a quadric and a form of degree $a$, for any $a$.
Let us begin with curves of degree 2.

\begin{ex}
The general plane section of any reduced curve $C$ of degree 2 is a reduced degree~2 zero-dimensional 
scheme, hence a complete intersection. If $C$ is connected, then it is a plane curve, hence aCM.
If $C$ is disconnected, then it consists of two skew lines, so it's non-aCM.

We observe that in this case, assuming that the curve is connected ensures its Cohen-Macaulayness.
\end{ex}

The situation is different for curves of degree $2a$, for $a\geq 2$.

\begin{ex}\label{rat}
Consider a (general) smooth rational curve $C$ of degree $2a$, $2\leq a$, lying on a smooth quadric 
surface $\Q\subset\PP^3$, e.g. the curve of parametric equations
$$\left\{ \begin{array}{l}  
x_0=s^{2a} \\
x_1=s^{2a-1}t \\
x_2=st^{2a-1} \\
x_3=t^{2a}
\end{array} \right.$$

$C$ is a rational, non-degenerate, smooth curve lying on the smooth quadric surface 
$\Q=x_0x_3-x_1x_2$. In fact, the saturated ideal of $C$ is
$$I_C=(x_0x_3-x_1x_2,x_0^{2a-2}x_2-x_1^{2a-1},x_0^{2a-3}x_2^2-x_1^{2a-2}x_3,\ldots,
x_0x_2^{2a-2}-x_1^2x_3^{2a-3},x_2^{2a-1}-x_1x_3^{2a-2}).$$
$C$ is non-aCM, since it has genus $g=0$, hence some entry of the $h$-vector has to be negative. 
In fact the only aCM, smooth rational curve in $\PP^3$  is the twisted cubic (general rational 
curve of degree~3).

Let $X$ be the general plane section of $C$. $X$ lies on a smooth conic and its $h$-polynomial is 
$h(z)=1+2z+2z^2+\ldots+2z^{a-1}+z^a$, since $X$ has the Uniform Position Property (see \cite{H}, about the 
$h$-vector of points in the plane with the UPP). Then $X$ is a complete intersection of type $(2,a)$.
\end{ex}

\begin{rem}\label{rat2}
In some cases, it will be useful to consider rational smooth curves, whose ideal is generated in small
degree. If $a$ is even, consider the curve $C$ of parametric equations
$$\left\{ \begin{array}{l}  
x_0=s^{2a} \\
x_1=s^{a+1}t^{a-1} \\
x_2=s^{a-1}t^{a+1} \\
x_3=t^{2a}
\end{array} \right.$$
Its saturated ideal is 
$$I_C=(x_0x_3-x_1x_2,x_0^2x_2^{a-1}-x_1^{a+1},x_0x_2^a-x_1^ax_3,x_2^{a+1}-x_1^{a-1}x_3^2).$$
If $a$ is odd, let $C$ be the curve parametrized by
$$\left\{ \begin{array}{l} 
x_0=s^{2a} \\
x_1=s^{a+2}t^{a-2} \\
x_2=s^{a-2}t^{a+2} \\
x_3=t^{2a}
\end{array} \right.$$
whose saturated ideal is 
$$I_C=(x_0x_3-x_1x_2,x_0^4x_2^{a-2}-x_1^{a+2},x_0^3x_2^{a-1}-x_1^{a+1}x_3,\ldots,
x_2^{a+2}-x_1^{a-2}x_3^4).$$
In both cases, $C$ is a rational, non-degenerate, smooth curve lying on the smooth quadric surface
$\Q=x_0x_3-x_1x_2$. As in Example~\ref{rat}, $C$ is non-aCM and its general plane section is a 
$CI(2,a)$.
The ideal of $I_C$ is generated in degree less than or equal to $a+1$ if $a$ is
even, and less than or equal to $a+2$ if $a$ is odd.
\end{rem}

The result of Strano has been generalized to curves in $\PP^{n+1}$ by R. Re in \cite{R}. It has been
further generalized to curves in $\PP^{n+1}$ with Gorenstein general hyperplane section, by C. Huneke and
B. Ulrich (see \cite{HU}). We will discuss their result extensively in the following section.
In \cite{M2}, J. Migliore proved a further generalization of their result for the case of hypersurface
section of a curve $C\subset\PP^{n+1}$.

\section{Conditions for Cohen-Macaulayness of a space curve}

In this section, we will be interested in finding conditions on the general hyperplane section 
of a curve $C$, that ensure Cohen-Macaulayness of the curve.
For the case when $X$ is a complete intersection, we refer to \cite{S} and \cite {R}.

We work over an algebraically closed field $k$ of characteristic $0$. The characteristic~$0$ hypothesis
is needed in Theorem~\ref{socle} of Huneke and Ulrich, and in its applications (Corollaries~\ref{bigger3}
and \ref{bigger3P^n}). Every other result and construction in this section is true over an algebraically
closed field of arbitrary characteristic.

In $\PP^2$, every arithmetically Gorenstein zero-dimensional scheme is a complete intersection. This is not 
the case in higher codimension, i.e. for zero-dimensional schemes in $\PP^n$ when $n\geq 3$.
The problem of finding a sufficient condition for a curve in $\PP^n$ to be arithmetically 
Cohen-Macaulay, hence arithmetically Gorenstein, given that its general hyperplane section is
arithmetically Gorenstein, has been solved by C.~Huneke and B.~Ulrich in \cite{HU}. This remarkable paper 
is based on a Lemma called the Socle Lemma; the Theorem that follows is a consequence of it, and we will 
make a substantial use of it in the sequel.

\begin{thm}\label{socle} (Theorem 3.16, \cite{HU})

Let $S=k[x_0,\ldots,x_{n+1}]$, $k$ a field of characteristic 0. Let $J\subset S$ be the homogeneous
ideal  of a reduced, connected curve $C\subset\PP^{n+1}$.
Let $L$ be a general linear form in $S$ and $X$ be the corresponding general hyperplane section of
$C$, $X\subset\PP^n$. The homogeneous ideal of $X$ in $R=S/(L)$ is $I=H^0_*(J+(L)/(L))\supseteq J+(L)/(L)$. Let 
$$0\lra \bigoplus_{i=1}^{b_{n-1}}R(-m_{n-1,i})\lra \cdots \lra \bigoplus_{i=1}^{b_1}R(-m_{1,i})\lra
I\lra 0$$ be the minimal free resolution of $I$ as an $R$-module. If $I\not= J+(L)/(L)$, then 
$$\mmin \{m_{n-1,i} \} \leq b+n-1$$ where $b=\mmin \{d\; |\; I_d\not= (J+(L)/(L))_d \}$.
\end{thm}

\begin{rem}
The curve $C$ is aCM if and only if $I= J+(L)/(L)$.

If $C$ is non-aCM, then there exists a minimal generator of $I=I_X$ of degree $b$, that is not the
image of any element of $J=I_C$ under the standard projection to the quotient. 
\end{rem}

\begin{rem}
It was observed by J. Migliore (see \cite{M2}, Proposition 2.2 and Theorem 2.4) that the hypotheses that
the curve
$C$ is reduced and connected are not necessary. In fact, one can show Theorem~\ref{socle} for any
curve $C\subset\PP^{n+1}$ that is non-degenerate, locally Cohen-Macaulay and equidimensional.

Notice moreover that the hypothesis on $C$ cannot be weakened any further. In fact, any
non-equidimensional curve is automatically non-aCM. Moreover, the general hyperplane section of a curve
only depends on its one-dimensional components. The hypothesis that $C$ is locally Cohen-Macaulay is
equivalent to $\M_C$ being of finite length as an $S$-module.
\end{rem}

Let us fix some notation. We will start with an analysis of the case of space curves.

Let $C\subset \PP^3$ be a curve, let $X\subset\PP^2$ be its general plane section.
Let $A$ be the homogeneous matrix whose maximal minors generate $I_X$ and $M$ be its degree matrix.
The minimal free resolution of $X$ is
$$0 \lra \bigoplus_{i=1}^t R(-m_i) \stackrel{A'}{\lra} \bigoplus_{j=1}^{t+1} R(-d_j) \lra I_X \lra 0$$
where  $d_1\geq d_2\geq\ldots\geq d_{t+1}$
are the degrees of a minimal system of generators, $m_1\geq m_2\geq\ldots\geq m_t$ and $A'$ is the
transpose of $A$.

The result that follows has been observed by J. Migliore in \cite{M2} (Proposition 2.2 and Remark 2.3),
and is an easy consequence of Theorem~\ref{socle}. 

\begin{cor}\label{bigger3}
Let $C\subset \PP^3$ be a curve, whose general plane section $X\subset\PP^2$ has
degree matrix
$M=(a_{i,j})_{i=1,\ldots ,t;j=1,\ldots ,t+1}$. If $a_{t,1}\geq 3$, then $C$ is arithmetically 
Cohen-Macaulay.
\end{cor}

\begin{proof} 
Let $L$ be the equation of the plane of $\PP^3$ in which $X$ is contained. $L$ 
is unique by non-degeneracy of $X$ and $C$. Assume by contradiction that $C$ is non arithmetically 
Cohen-Macaulay and let $b$ be the minimum degree in which the
ideal $I_X\subset S/(L)$ differs from $I_C+(L)/(L)\subset S/(L)$, as in the statement of
Theorem~\ref{socle}. By Theorem~\ref{socle} we have that
$$b\geq min \{m_i\}-2=m_t-2=d_1+a_{t,1}-2\geq d_1+1.$$
Hence all the minimal generators of $I_X$ come from images of the minimal 
generators of $I_C$. Then $C$ is arithmetically Cohen-Macaulay, contradicting our assumption.
\end{proof}

We will show in the sequel that the condition $a_{t,1}\geq 3$ is 
optimal. In fact, in Theorem~\ref{main} and Theorem~\ref{mmain} we will construct a reduced, connected,
non-aCM curve $C$ whose general plane section has degree matrix $M$, for any matrix $M$ with
$a_{t,1}\leq 2$.

We start with an analysis of the degree matrices corresponding to generic points.

\begin{ex}\label{mat1} (Degree matrix of three generic points)

Consider the degree matrix
$$M=\left(\begin{array}{ccc} 1 & 1 & 1 \\ 1 & 1 & 1 \end{array}\right).$$
$M$ is the degree matrix of three generic points in $\PP^2$.
A connected, reduced cubic curve $C\subset\PP^3$ is arithmetically Cohen-Macaulay. In fact, up to
isomorphism, the only integral, non-degenerate cubic curve in $\PP^3$ is the twisted cubic, that is aCM.
Any reducible, connected cubic curve is the union of a line and a plane conic (possibly reducible), meeting
in a point. The curves cannot lie on the same plane, otherwise the points of a general section of $C$ 
would be collinear.
Each of these curves is aCM. So it is not  possible to find a connected, reduced, non-aCM curve
$C\subset\PP^3$, whose  general plane section has degree matrix $M$. 

Dropping the requirement that the curve is connected, we can take $C$ to be the union of three skew 
lines in $\PP^3$, or the generic union of a line and a plane conic.
$C$ is smooth, disconnected and not arithmetically Cohen-Macaulay.

We also have a non-reduced curve: a fat line, whose ideal is given by $(L_1,L_2)^2$,
where $L_1, L_2$ are linearly independent linear forms. A fat line is a degree 3, non-degenerate aCM
curve. Its general plane section is a fat point, whose degree matrix is $M$.

For this particular matrix $M$ then, requiring that $C$ is connected forces Cohen-Macaulayness of 
the curve. 
Notice that Cohen-Macaulayness in this case does not follow from
Theorem~\ref{socle}.
\end{ex}

\begin{ex}(Generic points)\label{genpoints}

Let $X$ consist of $d$ generic points in $\PP^3$. The $h$-vector of $X$ is
$$h(z)=1+2z+\ldots+nz^{n-1}+\left( d-
\bin{n+1}{2}\right) z^n$$
where $n=\mmax\;\{\;i\; |\; \bin{i+1}{2}\leq d\;\}$. Let $s=d-\bin{n+1}{2}$. 
The initial degree of the saturated ideal 
$I_X$ is $\alpha (I_X)=n$, and the minimal free resolution of $I_X$ is
$$0 \lra R(-n-2)^s\oplus R(-n-1)^{n-2s} \lra R(-n)^{n+1-s} \lra I_X \lra 0 \;\;\;\; \mbox{if}
\;\;  0\leq s\leq \left[ \frac{n}{2} \right] $$
or
$$0 \lra R(-n-2)^s \lra R(-n)^{n+1-s}\oplus R(-n-1)^{2s-n} \lra I_X \lra 0 \;\;\;\; \mbox{if}
\;\; 
\left[ \frac{n}{2} \right] \leq s\leq n,$$
where $\left[ \frac{n}{2} \right]=\mmax\{ m\in\ZZ \; | \; 2m\leq n\}$.
The degree matrix for $X$ is then 
$$M=\underbrace{\left(\begin{array}{ll}
\left.\begin{array}{ccccc} 2 & \cdots & \cdots & 2 \\
\vdots & & & \vdots \\
2 & \cdots & \cdots & 2 
\end{array}\right\} s \\
\left.\begin{array}{cccc}
1 & \cdots & \cdots & 1 \\
\vdots & & & \vdots \\
1 & \cdots & \cdots & 1 
\end{array}\right\} n-2s
\end{array}\right)}_{n+1-s}$$
or respectively
$$M=\begin{array}{ccc}
\underbrace{\left(\begin{array}{ccc} 
1 & \cdots & 1 \\ 
\vdots & & \vdots \\ 
\vdots & & \vdots \\
1 & \cdots & 1 
\end{array}\right.}_{2s-n}
&
\underbrace{\left.\begin{array}{ccc}
2 & \cdots & 2 \\
\vdots & & \vdots \\
\vdots & & \vdots \\
 2 & \cdots & 2 
\end{array}\right)}_{n+1-s} 
&
\left.
\begin{array}{c}
 \\ 
 \\
 \\ 
 \\
 \\ 
 \end{array}\right\}
\end{array} s$$

\textbf{Claim.} {\em The general plane section of a general rational (smooth) curve of $\PP^3$ of degree
$d$ is a generic  set of $d$ points in the plane.}

Let us consider a generic zero-dimensional scheme $X$ of degree $d$ in the plane.
We only need to consider the case $d\geq 4$, since for $d=1,2,3$ a general rational curve of degree $d$ is
respectively a line, a smooth plane conic and a twisted cubic. In all of those cases we know that the 
general plane section consists of generic points. 
Notice that for $d\leq 3$ a general rational (smooth) curve is arithmetically Cohen-Macaulay.

By a result of Ballico and Migliore (see \cite{BM}, Theorem 1.6), we know
that there exists a smooth rational curve of 
degree $d$ that has $X$ as a proper section. Then, a generic rational curve $C$ of the
same degree $d$ will have a generic zero-dimensional scheme of degree $d$ as its proper section. By
upper-semicontinuity, we can then conclude that a general hyperplane section of $C$ is a generic
zero-dimensional scheme of degree $d$.

For all the degree matrices $M$ that correspond to $d$ generic points in the plane, $d\geq 4$, 
we can find a smooth rational curve whose general plane section has degree matrix $M$.
A smooth, rational curve of degree $d$ and genus $g=0$, with $h$-vector $(1,h_1,\ldots,h_s)$, has
$0=g=h_2+2h_3+\ldots +(s-1)h_s$. Then it cannot be aCM unless $s=1$, since for an aCM curve $h_i\geq 0$
for all $i$. In this case, $C$ has degree $d=h_0+h_1\leq 3$. 
\end{ex}

We are now going to analyze the general case. We will start from matrices of size $2\times 3$ or, more
generally, matrices of any size with an assumption on one of the entries. See Example~\ref{mat1} for 
the necessity of the assumption that $M$ is not a $2\times 3$ matrix with all the entries equal to $1$.

\begin{thm}\label{main}
Let $M=(a_{i,j})$ be a degree matrix of size 
$t\times (t+1)$ such that $a_{r,r-1}\leq 2$, for some $r$. 
Assume $M$ is not a $2\times 3$ matrix with all  the entries equal to 1. Then there exists a reduced,
connected, non-aCM curve  $C\subset\PP^3$ whose general plane section $X\subset\PP^2$ has degree 
matrix $M$.
\end{thm}

\begin{proof} 
Consider the two submatrices of $M=(a_{i,j})_{i=1,\ldots,t;j=1,\ldots,
t+1}$, $$L_1=(a_{i,j})_{i=1,\ldots,r-1; j=1,\ldots,r-1} \;\;\;\;\; 
N=(a_{i,j})_{i=r,\ldots,t; j=r,\ldots,t+1}$$
where $r$ is an integer $2\leq r\leq t$, such that $a_{r,r-1}\leq 2$.
Let $$a=a_{1,1}+a_{r,t+1}+a_{r,r}+a_{r+1,r+1}+\ldots+a_{t,t}-a_{r,1}$$ and
let $L$ be the matrix obtained by adding to $L_1$ the column 
$$(a,a-a_{1,r-1}+a_{2,r-1},a-a_{1,r-1}+a_{3,r-1},\ldots,a-a_{1,r-1}+a_{r-1,r-1})^t$$ as 
the $r$-th column.

Notice that all the entries on the diagonal on $L$ are positive, since they coincide with the first $r-1$
entries of the diagonal of $M$. The entries on the diagonal of $M$ are positive by
Proposition~\ref{diag_entries}. 
Moreover, $a-a_{1,r-1}=a_{1,1}+a_{r,t+1}+a_{r,r}+a_{r+1,r+1}+\ldots+a_{t,t}-a_{r,1}-a_{1,r-1}=
a_{r,t+1}+a_{r,r}+a_{r+1,r+1}+\ldots+a_{t,t}-a_{r,r-1}\geq a_{r,r}+a_{r+1,r+1}+\ldots+a_{t,t}>0$, by
Proposition~\ref{diag_entries}. So $a>a_{1,r-1}$ and $L$ is  a degree matrix, with the convention on the
order of the entries that we put in the definition (entries decrease from top to bottom and from right to
left).
The entries on the diagonal of $N$ are positive as well, since they are a subset of the entries on the
diagonal of $M$. Then, both $L$ and $N$ are degree matrices.

Let us consider two reduced, connected, arithmetically Cohen-Macaulay curves $C_1,C_2\subset\PP^3$, 
with degree matrices $N,L$ respectively. Let $C_1,C_2$ be generic through a fixed (common) 
point $P$. We can assume that a generic curve with a prescribed degree matrix is reduced, by \cite{Ga} or
by Proposition 2.5 in \cite{BCG}. Moreover, we can assume that $C_1$ and $C_2$ are connected curves, since
for any degree matrix there is a connected, arithmetically Cohen-Macaulay curve associated to it (so, for
a given degree matrix $N$, we can take the curve $E$ to be the cone over the zero-dimensional scheme constructed as
in \cite{Ga} or in Proposition 2.5 of \cite{BCG}).  Under the assumption that the entries on the
subdiagonal of
$M$ are positive (i.e. if $a_{i+1,i}>0$ for all $i$), so are the entries on the subdiagonals of $L$ and
$N$. Then by a result of T. Sauer (see \cite{Sa}), we can assume that $C_1$ and $C_2$ are also smooth.

Let $C=C_1\cup C_2$ be the union of the two curves. $C$ is reduced, non-degenerate and 
connected by construction. It has two irreducible components, both of them smooth if the subdiagonal of
$M$ is positive. Moreover, the ideal $I_{C_1}+I_{C_2}$ is not saturated, since its saturation is the
homogeneous ideal of a point. Looking at the short exact sequence
$$0 \lra I_C \lra I_{C_1}\oplus I_{C_2} \lra I_{C_1}+I_{C_2} \lra 0$$
we have that $\M_C=H^0_*(\I_{C_1}+\I_{C_2})/(I_{C_1}+I_{C_2})\not= 0$, so $C$ is not
arithmetically Cohen-Macaulay.

Taking a general plane section of $C$, we get a zero-dimensional scheme $X\subset\PP^2$, with 
saturated 
homogeneous ideal $I_X$. As a scheme, $X=X_1\cup X_2$, where $X_1,X_2$ are general plane 
sections of $C_1,C_2$ respectively.
Let the minimal free resolutions of $X_1$ and $X_2$ be
$$0 \lra \FF_2 \lra \FF_1 \lra I_{X_1} \lra 0$$
and
$$0 \lra \GG_2 \lra \GG_1 \lra I_{X_2} \lra 0.$$
Let $F$ be a generator of minimal degree in a minimal system of generators of $I_{X_2}$, and
let $d=deg(F)=a_{1,1}+a_{2,2}+\ldots+a_{r-1,r-1}$ (notice that $d>0$ by Proposition~\ref{diag_entries}).
By generality of our choice of $C_1$ and $C_2$, we can assume that $F$ is non-zerodivisor
modulo $I_{X_1}$. Consider now the ideal $I_{X_1}+(F)$. It is an Artinian ideal of 
$R=k[x_0,x_1,x_2]$, with minimal free resolution
\begin{equation}\label{mfr}
0 \lra \FF_2(-d) \lra \FF_2\oplus\FF_1(-d) \lra \FF_1\oplus R(-d) \lra I_{X_1}+(F) \lra 0
\end{equation}
and socle in degree $s=a_{1,1}+a_{2,2}+\ldots+a_{t,t}+a_{r,t+1}-3$. 

All minimal generators of $I_{X_2}$, except for $F$, 
have degrees bigger or equal to
$$d-a_{1,r-1}+a=2a_{1,1}+a_{2,2}+\ldots+a_{r-1,r-1}-a_{1,r-1}+a_{r,t+1}+a_{r,r}+
a_{r+1,r+1}+\ldots+a_{t,t}-a_{r,1}=$$ $$=a_{1,1}+\ldots+a_{t,t}-a_{r,r-1}+a_{r,t+1}=s+3-a_{r,r-1}\geq
s+1,$$ by assumption that $a_{r,r-1}\leq 2$. 
Since $s$ is the socle degree of the quotient ring $R/I_{X_1}+(F)$, $$I_{X_1}+(F)=I_{X_1}+I_{X_2}.$$ 
Let $$0 \lra \HH_2 \lra \HH_1 \lra I_X \lra 0$$ be a minimal free resolution of $I_X$.
Applying the Mapping Cone construction to the short exact sequence
$$0 \lra I_X \lra I_{X_1}\oplus I_{X_2} \lra I_{X_1}+I_{X_2}=I_{X_1}+(F) \lra 0$$
yields the following free resolution for $I_{X_1}+(F)$:
\begin{equation}\label{free}
0 \lra \HH_2 \lra \HH_1\oplus\GG_2\oplus\FF_2 \lra \GG_1\oplus\FF_1 \lra I_{X_1}+(F) 
\lra 0.\end{equation}
Comparing (\ref{mfr}) and (\ref{free}) gives 
$$\HH_2=\GG_2\oplus\FF_2(-d)\oplus\FF, \;\;\;\;\; \HH_1=\GG_1'\oplus\FF_1(-d)\oplus\FF$$
for $\FF$ some free $R$-module and $\GG_1=\GG_1'\oplus R(-d)$. 
This follows from the fact that there can be no cancellation between $\GG'_1$ and $\FF_2$ in the
resolution of $I_{X_1}+(F)$
obtained via the Mapping Cone, since the two free modules come from the same minimal free resolution
(the one of 
$I_{X_1}\oplus I_{X_2}$). 
Moreover, the shifts of the free summand of $\GG_2$ are all different from the shifts 
of the free summands of $\FF_1(-d)$. In fact, the smallest shift among the free summands in $\GG_2$
is  $d+a+a_{r-1,r-1}-a_{1,r-1}=d+a_{t,r+1}+a_{r,r}+\ldots+a_{t,t}-a_{r,1}+a_{r-1,1}
>d+a_{t,r+1}+a_{r+1,r+1}+\ldots+a_{t,t}$, 
that is the highest shift among the free summands of $\FF_1(-d)$.

The free summands 
$\FF$ cannot split off in the minimal free resolution of $I_{X_1}+(F)$, because they
come from the minimal free resolution of $I_X$, hence the map between them is not an 
isomorphism on any free submodule (the map is left unchanged under the Mapping Cone). 
Then $\FF=0$, since (\ref{free}) has to equal (\ref{mfr}), after splitting. We obtain the following
minimal free resolution for $I_X$:
$$0 \lra \GG_2\oplus\FF_2(-d) \lra \GG_1'\oplus\FF_1(-d) \lra I_X \lra 0.$$
The degree matrix of $X$ is then $(b_{i,j})$, where $$b_{i,j}=a_{i,j}\;\; \mbox{for $1\leq i\leq r-1,
1\leq j\leq r-1$ and $r\leq i\leq t, r\leq j\leq t+1$.}$$ Moreover, $$b_{r,1}=d+\mbox{(maximum shift in
$\FF_2)-$(maximum shift in $\GG_1'$).}$$ Then 
$$b_{r,1}=d+(a_{r,r}+\ldots+a_{t,t}+a_{r,t+1})-(d-a_{1,1}+a)=a_{r,1}.$$
Notice that, since $M$ is homogeneous, all of its entries are determined by $L_1,N$ and $a_{r,1}$.
This proves that $M$ is the degree matrix of $X$.
\end{proof}

\begin{rem}\label{defmod}
We can easily compute the deficiency module of the curves constructed in Theorem~\ref{main}.
In fact, $C=C_1\cup C_2$ with $C_1$ and $C_2$ aCM meeting in exactly one point $P$. So we
have the exact sequence 
$$0 \lra I_C \lra I_{C_1}\oplus I_{C_2} \lra I_P \lra \M_C \lra \M_{C_1}\oplus\M_{C_2}=0$$
that together with the short exact sequence
$$0 \lra I_C \lra I_{C_1}\oplus I_{C_2} \lra I_{C_1}+I_{C_2} \lra 0$$
gives the isomorphism
$$\M_C\cong I_P/I_{C_1}+I_{C_2}.$$
In particular, $\alpha(\M_C)=1$.
\end{rem}

\begin{rem}\label{alt-const}
If, instead of taking $C_1$ and $C_2$ generic through the same point, we take them generic and disjoint 
(with the prescribed degree matrices), we get a non-degenerate, reduced, non-aCM, disconnected
curve with two connected components, 
$C=C_1\cup C_2$. A general plane section of $C$ has degree matrix $M$. The proof is very similar to that 
of Theorem~\ref{main}.

If the entries on the subdiagonal of $M$ are positive, we can take $C_1$ and $C_2$ to be smooth
and integral. In this case, $C$ is a non-degenerate, smooth, non-aCM, disconnected
curve with two smooth connected components, 
$C=C_1\cup C_2$, whose general plane section has degree matrix $M$.

In the case that $C_1$ and $C_2$ are disjoint, we can explicitly compute the deficiency
module $\M_C$ of $C$. We have the exact sequence 
$$0 \lra I_C \lra I_{C_1}\oplus I_{C_2} \lra R \lra \M_C \lra \M_{C_1}\oplus\M_{C_2}=0$$
that together with the short exact sequence
$$0 \lra I_C \lra I_{C_1}\oplus I_{C_2} \lra I_{C_1}+I_{C_2} \lra 0$$
gives the isomorphism
$$\M_C\cong R/I_{C_1}+I_{C_2}.$$ In particular $\alpha(\M_C)=0$.
\end{rem}

The construction of Theorem~\ref{main} is very simple in the case of matrices of size $2\times 3$. In
this case, moreover, $r=2$ and the condition $a_{r,r-1}=a_{2,1}\leq 2$ is always satisfied.

\begin{ex}\label{2x3case}
Consider a degree matrix 
$$M=\left(\begin{array}{ccc}
a & b & c \\
d & e & f \end{array}\right).$$
In order for $M$ to be a degree matrix, all the entries have to be positive, except possibly for $d$. By 
assumption $d\leq 2$.
Following the proof of Theorem~\ref{main}, let $C=CI(a,b+f)\cup CI(e,f)\subset\PP^3$, where the
complete intersections are generic through a common point $P$. Then $C$ is a non-aCM, connected, reduced,
non-degenerate space curve, smooth outside of $P$, whose general plane section has degree matrix $M$.
Moreover, the deficiency module is $\M_D\cong (x_1,x_2,x_3)/(F_1,F_2,G_1,G_2)$, where $(F_1,F_2)$ and
$(G_1,G_2)$ are the ideals of two generic complete intersections of type $(a,b+f)$ and $(e,f)$ through the
point $[1:0:0:0]$.

Let $D=CI(a,b+f)\cup CI(e,f)\subset\PP^3$, where the complete intersections are generic, hence disjoint. 
Then $D$ is a non-aCM, reduced space curve, with two smooth connected components. The general plane
section $X$ of $D$ has degree matrix $M$. Moreover, the deficiency module is $\M_D\cong
k[x_0,x_1,x_2,x_3]/CI(a,e,f,b+f)$.

Notice that in both cases, the initial degree of the ideal $I_C$ of $C$ is the same as the initial degree
of $I_X$, and the highest degree for a minimal generator of $I_C$ is $b+2f$.
\end{ex}

\begin{rem}
The assumption that $a_{r,r-1}\leq 2$ in Theorem~\ref{main} is essential. In fact, if $a_{r,r-1}>2$ 
then the size of the degree matrix that we obtain following the procedure of the Theorem is
strictly bigger than $t\times (t+1)$ (see Example~\ref{ar,r-1} below).

Nevertheless, any choice of $r\geq 2$ for which $a_{r,r-1}\leq 2$ would work. 
So, for a given degree matrix $M$, for any choice of $r$ such that $a_{r,r-1}\leq 2$ we
constructed a non-aCM curve, whose general plane section has degree matrix $M$. 
The curves that we get for two different values of $r$ are not projectively isomorphic, since their
connected components are not (in particular, their connected components have different degree matrices).
\end{rem}

In the next example, we show how the construction of Theorem~\ref{main} does not yield a curve whose
general plane section has the desired degree matrix, in the case that the hypothesis $a_{r,r-1}\leq 2$
is not  satisfied.

\begin{ex}\label{ar,r-1}
Let $$M=\left(\begin{array}{cccc} 
3 & 4 & 4 & 5 \\
2 & 3 & 3 & 4 \\
2 & 3 & 3 & 4
\end{array}\right)$$
and let $r=3$. Notice that $a_{3,2}=3\not\leq 2$. Let 
$$L=\left(\begin{array}{ccc} 
3 & 4 & 8 \\
2 & 3 & 7
\end{array}\right), \;\;\;\;\;\;\;\; N=(3,4)$$
and let $C_1,C_2$ be aCM, smooth, generic curves through a common point, with degree matrices $N, L$
respectively. Let $X_1, X_2$ be general plane sections of $C_1, C_2$, respectively.
The minimal free resolution of $I=I_{X_1}+I_{X_2}$, computed with \cite{CoCoA}, turns out to be  
$$\begin{array}{ccccc}
 & R(-7)\oplus R(-9) & & R(-3)\oplus R(-4) & \\
0\lra R(-12)^3 \lra & \oplus & \lra & \oplus & \lra I\lra 0 \\
 & R(-10)\oplus R(-11)^3 & & R(-6)\oplus R(-10) & \end{array}$$
hence the degree matrix of the general plane section of $C=C_1\cup C_2$ is
$$M'=\left(\begin{array}{cccccc}  3 & 3 & 3 & 3 & 4 & 5 \\
2 & 2 & 2 & 2 & 3 & 4 \\
1 & 1 & 1 & 1 & 2 & 3 \\
1 & 1 & 1 & 1 & 2 & 3 \\
1 & 1 & 1 & 1 & 2 & 3 
\end{array}\right)$$
and not the required matrix $M$.
The problem comes from the fact that the socle of $I_{X_1}+(F)$ has final degree $10$, and $I_{X_2}$ has
a minimal generator in degree $10$ that does not belong to $I_{X_1}+(F)$.
\end{ex}

\begin{rem}
In the statement of Theorem~\ref{main}, we pointed out that the 
construction does not work for the matrix 
$$M=\left(\begin{array}{ccc} 1 & 1 & 1 \\ 1 & 1 & 1 \end{array}\right)$$
that we analyzed in Example~\ref{mat1}.
In fact, for this matrix our construction yields curves $C_1=$ a plane conic and $C_2$= 
a line, meeting in a point. In this case, $I_{C_1}+I_{C_2}$ is saturated and $C=C_1\cup C_2$ 
is arithmetically Cohen Macaulay.

Notice that, if we take a generic (disjoint) union of a line and a conic, we get a non-degenerate,
smooth, non-aCM, disconnected curve, whose general plane section consists of three generic points and has
degree matrix $M$. 
\end{rem}

We now present an alternative construction for the degree matrices of size $2\times 3$. The advantage 
with respect to the construction of Theorem~\ref{main} is that the saturated ideal of the curves that we 
obtain in the following theorem are minimally generated in low degree. This will be useful in the next 
section.

\begin{thm}\label{2x3nice}
Let $M$ be a degree matrix of size $2\times 3$, 
$$M=\left(\begin{array}{ccc} 
a_{1,1} & a_{1,2} & a_{1,3} \\
a_{2,1} & a_{2,2} & a_{2,3} 
\end{array}\right)$$ and assume that $a_{2,1}\leq 2$. Then there exists a reduced,
connected, non arithmetically Cohen-Macaulay curve  $C\subset\PP^3$, whose general plane section 
$X\subset\PP^2$ has degree matrix $M$, and such that the saturated ideal $I_C$ of $C$ is minimally 
generated in degree smaller than or equal to $a_{1,2}+a_{2,3}+1$.
\end{thm}

\begin{proof}
Let $C_1$ be a generic complete intersection of type $(a_{2,2},a_{2,3})$, and let
$I_{C_1}\subset S=k[x_0,\ldots,x_3]$ be the saturated ideal of $C_1$. Let $G$ be a generic form of degree
$a_{1,1}$. Then $I=I_{C_1}+(G)$ is the saturated ideal of a generic complete intersection of type 
$(a_{1,1},a_{2,2},a_{2,3})$. Therefore the scheme $Z$ associated to $I$ is a zero-dimensional scheme, 
consisting of $a_{1,1}\cdot a_{2,2}\cdot a_{2,3}$ distinct points. Let $P$ be one of the points of $Z$, and let 
$X=Z-P$ be the complement of $P$ in $Z$. Notice that $X$ is linked to $P$ via the complete intersection 
$Z$, therefore by Proposition 5.2.10 in \cite{M} one 
gets a free resolution of $I_X$ of the form
$$\begin{array}{rcl}
 & S(-a_{1,1}-a_{2,2}-a_{2,3}+2)^3 \oplus & \\ 
0\lra S(-a_{1,1}-a_{2,2}-a_{2,3}+1)^3\lra & S(-a_{1,1}-a_{2,2})\oplus & \lra \\
 & S(-a_{2,2}-a_{2,3})\oplus S(-a_{1,1}-a_{2,3}) &
\end{array}$$
$$\begin{array}{rcl} 
 & S(-a_{1,1}-a_{2,2}-a_{2,3}+3)\oplus & \\
\lra & S(-a_{1,1})\oplus S(-a_{2,2})\oplus & \lra I_X\lra 0.\\
 & S(-a_{2,3}) &
\end{array}$$
The resolution is not a priori minimal.

The socle of the complete intersection $Z$ is concentrated in degree\newline 
$a_{1,1}+a_{2,2}+a_{2,3}-3\leq a_{1,2}
+a_{2,3}-1$, since $a_{2,1}\leq 2$ by assumption. Therefore, the Hilbert function of $Z$ in degree 
$a_{1,2}+a_{2,3}$ is $$H_Z(a_{1,2}+a_{2,3})=deg(Z).$$ The Hilbert function of $X$ in the same degree is
$$H_X(a_{1,2}+a_{2,3})\leq deg(X)=deg(Z)-1.$$ Then there is a surface $F$ of degree $a_{1,2}+a_{2,3}$ 
that contains $X$ but does not contain $Z$, so it contains $X$ and not $P$. Let the surface $F$ be generic, 
with this property.
Let the curve $C_2$ be the scheme-theoretic intersection of $F$ and $G$. 
$C_2$ is a complete intersection curve of type $(a_{1,1},a_{1,2}+a_{2,3})$.
By the construction, $C_1\cap C_2=X$. Let $C$ be the union of the two complete intersection curves, 
$C=C_1\cup C_2$. $C$ is reduced and connected, and it has two irreducible components. 
Its general plane section is the union of 
a $CI(a_{1,1},a_{1,2}+a_{2,3})$ and a $CI(a_{2,2},a_{2,3})$. The same argument as in the proof of 
Theorem~\ref{main} applies, showing that the general plane section of $C$ has degree matrix $M$.

We need to show that $C$ is not arithmetically Cohen-Macaulay. 
From the long exact sequence 
$$0\lra I_C\lra I_{C_1}\oplus I_{C_2}\lra I_X\lra M_C\lra 0$$ we see that the deficiency module of $C$ is 
$$\M_C\cong I_X/(I_{C_1}+I_{C_2}).$$ Then $C$ is arithmetically Cohen-Macaulay if and only if 
$I_{C_1}+I_{C_2}=I_X$, if and only if $I_{C_1}+I_{C_2}$ is saturated.

In order to show that the ideal $I_{C_1}+I_{C_2}$ is not saturated, we compute a free resolution of it. 
Multiplication by $F$ in $S/I$ yields the long exact sequence
$$0\lra (I:F)/I(-a_{1,2}-a_{2,3})\lra S/I(-a_{1,2}-a_{2,3})\lra S/I\lra S/(I+(F))\lra 0.$$
$I:F=I:(I+(F))$, and since $I+(F)=I_{C_1}+I_{C_2}$, then $I:F=I:(I_{C_1}+I_{C_2})$.
The saturation of $I_{C_1}+I_{C_2}$ is $I_X$, since $C_1\cap C_2=X$.
Therefore $$I:F=I:(I_{C_1}+I_{C_2})=I:I_X=I_P.$$ 
The last equality follows from the fact that $P$ is the residual 
to $X$ in the complete intersection $Z$, whose homogeneous saturated ideal is $I$.
Then $$I:F=I_P\;\;\;\;\;\mbox{and}\;\;\;\;\; I_{C_1}+I_{C_2}=I+(F).$$ These equalities give 
the short exact sequence 
$$0\lra S/I_P(-a_{1,2}-a_{2,3})\lra S/I \lra S/(I_{C_1}+I_{C_2})\lra 0.$$ Using the Mapping Cone 
construction, we obtain the free resolution for $I_{C_1}+I_{C_2}$
$$\begin{array}{rccl}
 & S(-a_{1,2}-a_{2,3}-2)^3 & & \\
0\lra S(-a_{1,2}-a_{2,3}-3) \lra & \oplus & \lra & \\
 & S(-a_{1,1}-a_{2,2}-a_{2,3}) & & 
\end{array}$$
$$\begin{array}{rcccl}
 & S(-a_{1,2}-a_{2,3}-1)^3\oplus & & S(-a_{1,1})\oplus & \\
\lra & S(-a_{1,1}-a_{2,2})\oplus & \lra & S(-a_{2,2})\oplus & \lra I_{C_1}+I_{C_2}\lra 0.\\
 & S(-a_{2,2}-a_{2,3})\oplus & & S(-a_{2,3})\oplus & \\
 & S(-a_{1,1}-a_{2,3}) & & S(-a_{1,2}-a_{2,3}) & 
\end{array}$$
The resolution is not minimal a priori. However, no cancellation can take place between the last free module and 
the following one, because $a_{1,2}+a_{2,3}+3>a_{1,1}+a_{2,2}+a_{2,3}$, since $a_{2,1}<3$. This proves that 
the ideal $I_{C_1}+I_{C_2}$ is not saturated, and therefore $C$ is not arithmetically Cohen-Macaulay.

Consider the short exact sequence 
$$0\lra I_{C_1}+I_{C_2} \lra I_X \lra \M_C\lra 0.$$
The Mapping Cone procedure produces a free resolution of $\M_C$ of the form 
$$\begin{array}{rccl}
 & S(-a_{1,2}-a_{2,3}-2)^3 & & \\
0\lra S(-a_{1,2}-a_{2,3}-3) \lra & \oplus & \lra & \\
 & S(-a_{1,1}-a_{2,2}-a_{2,3}) & & 
\end{array}$$
$$\begin{array}{rcccl}
 & S(-a_{1,2}-a_{2,3}-1)^3\oplus & & S(-a_{1,2}-a_{2,3})\oplus S(-a_{1,1})\oplus S(-a_{2,2}) & \\ 
\ra & S(-a_{1,1}-a_{2,2})\oplus S(-a_{1,1}-a_{2,3}) & \ra & S(-a_{2,2}-a_{2,3})\oplus S(-a_{1,1}-a_{2,3}) & \ra \\
 & S(-a_{2,2}-a_{2,3}) & & S(-a_{2,3})\oplus S(-a_{1,1}-a_{2,2}) & \\
 & S(-a_{1,1}-a_{2,2}-a_{2,3}+1)^3 & & S(-a_{1,1}-a_{2,2}-a_{2,3}+2)^3 &
\end{array}$$
$$\begin{array}{rcl}
 & S(-a_{1,1}-a_{2,2}-a_{2,3}+3)\oplus & \\
\lra & S(-a_{2,3})\oplus &\lra \M_C\lra 0.\\
 & S(-a_{1,1})\oplus S(-a_{2,2}) & 
\end{array}$$
The free summands $S(-a_{1,1})\oplus S(-a_{2,2})\oplus S(-a_{2,3})$ in the first free module of the 
resolution of $\M_C$ come from the free resolution of $I_X$. Since the minimal generators of $I_{C_1}+I_{C_2}$ 
in those degrees coincide with the minimal generators of $I_X$, the free summands that did not already 
cancel in the minimal free resolution of $I_X$ cancel in the minimal free resolution of $\M_C$ with the 
corresponding free summands in the second free module (coming from the free resolution of 
$I_{C_1}+I_{C_2}$). Therefore the first free module in the minimal free resolution of $\M_C$ is simply 
$S(-a_{1,1}-a_{2,2}-a_{2,3}+3)$. This proves that the initial degree of $\M_C$ is 
$$\alpha(M_C)=a_{1,1}+a_{2,2}+a_{2,3}-3.$$
From the shifts in the free resolution of $\M_C$, one can also deduce an upper bound for the 
Castelnuovo-Mumford regularity of $\M_C$:
$$reg(\M_C)\leq a_{1,2}+a_{2,3}-1.$$
From Lemma 3.12 in \cite{GM} it follows that, since the saturated ideal of the general hyperplane section of 
$C$ has no minimal generators in degree bigger than or equal to $a_{1,2}+a_{2,3}+1$, and the last non-zero 
component of the deficiency module of $C$ occurs in degree $$\alpha^+(\M_C)\leq a_{1,2}+a_{2,3}-1,$$
then the ideal $I_C$ is minimally generated in degree smaller than or equal to $a_{1,2}+a_{2,3}+1$.
\end{proof}

\begin{rem}\label{MCnice}
In the proof of Theorem~\ref{2x3nice} we compute a free resolution of the deficiency module $\M_C$ of the 
curve $C$ that we construct. Moreover, we prove that
$$\alpha(M_C)=a_{1,1}+a_{2,2}+a_{2,3}-3, \;\;\;\;\; \alpha^+(\M_C)\leq a_{1,2}+a_{2,3}-1$$ and that the 
Castelnuovo-Mumford regularity of $\M_C$ is bounded by $$reg(\M_C)\leq a_{1,2}+a_{2,3}-1.$$
\end{rem}

\begin{rem}
The saturated ideal of the general plane section $X$ of the curve $C$ has a minimal generator in 
degree $a_{1,2}+a_{2,3}$. Therefore the ideal of any curve that has $X$ as a general plane section 
necessarily has a minimal generator in degree $a_{1,2}+a_{2,3}$ or higher. 
\end{rem}

For the arguments that follow, we need to show the existence of smooth surfaces containing the curves 
constructed in Theorem~\ref{2x3nice}.

\begin{lemma}\label{smnice}
Let $C$ be a curve constructed as in Theorem~\ref{2x3nice}. For each $d\geq a_{1,2}+a_{2,3}+1$ there is a 
smooth surface of degree $d$ containing $C$. 
\end{lemma}

\begin{proof}
Consider the linear system $\Delta$ of surfaces of $\PP^3$ of degree $d$ containing $C$.
$C=C_1\cup C_2$ is a union of $2$ complete intersection curves. 
Let $Sing(C)=X\cup Y$ be the singular locus of $C$. $C$ is singular at the points where the two components
intersect, and possibly at some other zero-dimensional subset $Y\subset C_2$. 
The general element of $\Delta$ is basepoint-free outside of $C$, hence smooth
outside of $C$ by Bertini's Theorem. Now consider a point $P\in C$. We want to show that the 
general element of $\Delta$ is smooth at $P$. By Corollary 2.10 in \cite{GV},
it is enough to exhibit two elements of $\Delta$ meeting transversally at $P$. Since $C$
is smooth outside of $Sing(C)$, for each point $P\not\in Sing(C)$ we have two minimal generators of
$I_C$, call them $F$ and $G$, meeting transversally at $P$. The degree of each of them is at most $d$. 
Add generic planes as needed, to obtain surfaces of degree $d$ that meet transversally at $P$. 

In order to complete the proof, we need to check that the points of $Sing(C)$ are not fixed singular
points for $\Delta$. So it is enough to find a surface for each $P\in Sing(C)$ that contains $C$ and is
non-singular at $P$. For each point $Q\in Y$ we have a smooth surface $G$ containing $C_2$. Taking the union
of $G$ with a smooth surface of $C_1$ of appropriate degree ($C_1$ is smooth, so we can always find 
such a surface) that does not contain $Q\not\in C_1$ gives a surface that is smooth at $Q$ and contains $C$.

Let $Q\in X$. We need to find a surface containing $C$ that is smooth at $Q$.
Let $F_1,\ldots,F_n$ be a minimal system of generators of $I_C$. $d_i:=deg(F_i)\leq d$ for all $i$.
Some of the minimal generators 
of $I_C$ are smooth at $Q$ (the ones of degrees $a_{1,1},a_{2,2},a_{2,3}$ are smooth by genericity). Assume 
that $F_1,\ldots,F_r$ are smooth at $Q$. 
Then let $T=G_1F_1+\ldots+G_rF_r$ where each $G_i$ is a generic polynomial of degree $d-d_i$.
The surface defined by $T$ contains $C$ by construction. In order to check that $T$ is smooth at $Q$, it 
suffices to show that not all the partial derivatives of $T$ vanish at $Q$. 
Denote the derivative of $F_i$ with respect to $x_j$ by $F_{i,j}$.
Some of the partial derivatives of $F_i$ do not vanish at $Q$.
For example, assume that $F_{1,2}(Q)\neq 0$. Then the partial derivative of $T$ with respect to $x_2$ 
evaluated at $Q$ is $T_2(Q)=
G_{1,2}(Q)F_1(Q)+\ldots+G_{r,2}(Q)F_r(Q)+G_1(Q)F_{1,2}(Q)+\ldots+G_r(Q)F_{r,2}(Q)=
G_1(Q)F_{1,2}(Q)+\ldots+G_r(Q)F_{r,2}(Q).$ By genericity of $G_1,\ldots,G_r$ we can assume that none 
of them vanishes at $Q$ and that $G_1(Q)F_{1,2}(Q)+\ldots+G_r(Q)F_{r,2}(Q)\neq 0$.
This shows smoothness of $T$ at $Q$, and therefore concludes the proof.
\end{proof}

For any degree matrix $M$, such that one of its entries is smaller than or equal to $2$, 
we are going to construct an example of a reduced, connected, non-aCM curve, whose general plane section
has degree  matrix $M$. Notice that not all degree matrices can correspond to points that are the general
plane section of an integral curve. In particular, none of the curves that we will construct in the proof
of the Theorem will be integral. We will deal with degree matrices of points that can lift to an integral
curve in the next section.

\begin{thm}\label{mmain}
Let $M=(a_{i,j})$ be a degree matrix of size $t\times (t+1)$ such that $a_{t,1}\leq 2$. 
Assume $M$ is not a $2\times 3$ matrix with all  the entries equal to 1. Then there exists a reduced,
connected, non-aCM curve $C\subset\PP^3$ whose general plane section $X\subset\PP^2$ has degree 
matrix $M$.
\end{thm}

\begin{proof}
We will proceed by induction on the size $t$ of $M$. We will include in the induction hypothesis that
$\alpha(I_C)=\alpha(I_X)$. The thesis is true for $t=2$, as shown in Theorem~\ref{main} and in
Example~\ref{2x3case}. In fact, from the proof of the Theorem it follows that
$\alpha(I_C)=a_{1,1}+a_{2,2}$, since
$C=CI(a_{1,1},a_{1,2}+a_{2,3})
\cup CI(a_{2,2},a_{2,3})$.

Let $M=(a_{i,j})_{i=1,\ldots t;\; j=1,\ldots t+1}$ be a degree matrix with $a_{t,1}\leq 2$. Assume that
$a_{t-1,1}\leq 2$ and let $N=(a_{i,j})_{i=1,\ldots t-1;\; j=1,\ldots t}$ be the submatrix of $M$
consisting of the first 
$t-1$ rows and the first $t$ columns. $N$ is a degree matrix, since the entries on its diagonal
agree with the first $t-1$ entries on the diagonal of $M$, so they are positive. By the induction 
hypothesis, we have a non-aCM, reduced, connected curve $D\subset\PP^3$, whose general plane section
$Y\subset\PP^2$ has degree matrix $N$. Moreover, $\alpha(I_D)=\alpha(I_Y)=a_{1,1}+\ldots, a_{t-1,t-1}$.
Let $S$ be a surface of degree $s=a_{1,1}+\ldots+a_{t,t}$ containing $D$. Such an $S$ exists,
since
$s=\alpha(I_D)+a_{t,t}>\alpha(I_D)$. Moreover, $S$ can be chosen such that its image in $I_Y$ is not a minimal
generator, since $deg(S)>\alpha(I_Y)$. Perform a basic double link on
$S$, with a generic surface
$F$ of degree $a_{t,t+1}>0$, that meets $D$ in (at least) a point. Let $C=D\cup(S\cap F)$. Then $C$ is
reduced and connected, and $\M_C\cong \M_D(-a_{t,t+1})\neq 0$, so $C$ is non-aCM. $C$ is non-degenerate,
since $D$ is not. Moreover, by generality of our choices, $D$ and $S\cap F$ meet transversally at each
of their points of intersection, and each of their points of intersection is a smooth point on both $D$
and $S\cap F$. We have the short exact sequence (see
\cite{M}, Proposition 5.4.5)
$$0\lra R(-s-a_{t,t+1})\lra I_Y(-a_{t,t+1})\oplus R(-s)\lra I_X \lra 0.$$
Then, using Mapping Cone, a free resolution of $I_X$ is given by
$$\begin{array}{ccccc}
 & R(-s-a_{t,t+1}) & & R(-s) & \\
0\lra & \oplus & \lra & \oplus & \lra I_X\lra 0. \\
 & \FF_2(-a_{t,t+1}) & & \FF_1(-a_{t,t+1}) \end{array}$$
Notice that $R(-s-a_{t,t+1})$ cannot split off with any of the free summands of $\FF_1(-a_{t,t+1})$, since
$S$ is not a minimal generator. Moreover, none of the shifts appearing in $\FF_2(-a_{t,t+1})$
can be equal to $s$, since $a_{i,t+1}>0$ for all $i$ (if $a_{i,t+1}\leq 0$ for some $i$, then $a_{i,j}\leq
0$ for all $j$, and this is not possible for a degree matrix). This shows that the resolution is
minimal. Then, the degree matrix of $X$ is $M$, as required. Moreover, $\alpha(I_C)\leq s=\alpha(I_X)$,
so $\alpha(I_C)=\alpha(I_X)$.

The case when $a_{t-1,1}\geq 3$ and $a_{t,t-1}>0$ is analogous: let $N=(a_{i,j})_{i=2,\ldots t;\; 
j=1,\ldots t}$ be the submatrix of $M$ consisting of the last 
$t-1$ rows and the first $t$ columns. Notice that since $a_{t-1,1}\geq 3$, then $a_{i+1,i}>0$ for
$i=1,\ldots,t-2$. Perform a basic double link on a surface $S$ of degree $a_{1,1}+\ldots+a_{t,t}$ with a
form $F$ of degree $a_{1,t+1}$ (see \cite{M} about basic double links).

The case when $a_{t-1,1}\geq 3$ and $a_{t,t-1}\leq 0$ is again similar. Let $N=(a_{i,j})_{i=2,\ldots
t;\;  j=2,\ldots t+1}$ be the submatrix of $M$ consisting of the last $t-1$ rows and the last $t$
columns. Notice that $a_{t,2}\leq 0<2$. Perform a basic double link on a surface $S$ of degree
$a_{1,2}+\ldots+a_{t,t+1}$ with a form $F$ of degree $a_{1,1}$.
\end{proof}

\begin{rem}\label{components}
The curve $C$ that we constructed in Theorem~\ref{mmain} is a union of $t$ complete intersections. More
precisely, if $a_{k,l}\leq 2$ and $a_{k-1,l}>0$, $a_{k,l+1}>0$, then 
$C$ can be built following the inductive procedure we showed, starting from the submatrix
$$\left( \begin{array}{ccc} 
a_{k-1,l} & a_{k-1,l+1} & a_{k-1,l+2} \\
a_{k,l} & a_{k,l+1} & a_{k,l+2}
\end{array} \right).$$
Notice that one can always find such $k,l$. Moreover, one can assume that $l\leq k-1$, since the entries 
on the diagonal of $M$ are positive.

Then $C$ is the union
$$C=CI(a_{k-1,l},a_{k-1,l+1}+a_{k,l+2})\cup CI(a_{k,l+1},a_{k,l+2})\cup$$
$$CI(a_{k-2,l}+a_{k-1,l+1}+a_{k,l+2},a_{k-2,l-1})\cup\ldots\cup
CI(a_{k-l,2}+\ldots+a_{k,l+2},a_{k-l,1})\cup$$
$$CI(a_{k-l-1,1}+\ldots+a_{k,l+2},a_{k-l-1,l+3})\cup
\ldots\cup CI(a_{1,1}+\ldots+a_{k,k},a_{1,k+1})\cup$$
$$CI(a_{1,1}+\ldots+a_{k+1,k+1},a_{k+1,k+2})\cup\ldots\cup CI(a_{1,1}+\ldots+a_{t,t},a_{t,t+1}).$$

If $l\leq k-2$, then $C$ can be taken to be the union
$$C=CI(a_{k-1,l},a_{k-1,l+1}+a_{k,l+2})\cup CI(a_{k,l+1},a_{k,l+2})\cup$$
$$CI(a_{k-2,l}+a_{k-1,l+1}+a_{k,l+2},a_{k-2,l-1})\cup\ldots\cup
CI(a_{k-l,2}+\ldots+a_{k,l+2},a_{k-l,1})\cup$$
$$CI(a_{k-l-1,1}+\ldots+a_{k,l+2},a_{k-l-1,l+3})\cup
\ldots\cup CI(a_{1,1}+\ldots+a_{k,k},a_{1,k+1})\cup$$
$$CI(a_{1,1}+\ldots+a_{k+1,k+1},a_{k+1,k+2})\cup\ldots\cup CI(a_{1,1}+\ldots+a_{t,t},a_{t,t+1}).$$

If $l=k-1$, then $C$ can be taken to be the union
$$C=CI(a_{k-1,k-1},a_{k-1,k}+a_{k,k+1})\cup CI(a_{k,k},a_{k,k+1})\cup$$
$$CI(a_{k-2,k-1}+a_{k-1,k}+a_{k,k+1},a_{k-2,k-2})\cup\ldots\cup
CI(a_{1,2}+\ldots+a_{k,k+1},a_{1,1})\cup$$
$$CI(a_{1,1}+\ldots+a_{k+1,k+1},a_{k+1,k+2})\cup\ldots\cup CI(a_{1,1}+\ldots+a_{t,t},a_{t,t+1}).$$

Clearly there are other ways to perform the basic double links other than the examples that we present here. 
Different sequences of basic double links yield curves that are not projectively isomorphic, since they are 
unions of complete intersections of different types. Therefore, following the construction of 
Theorem~\ref{mmain}, one can produce different curves from the examples that we just gave. 

One can easily show by induction that $$\M_C\cong (L_1,L_2,L_3)/(F_1,F_2,G_1,G_2)(-a)$$ 
as an $S$-module, where $$a=a_{k-2,l-1}+\ldots+a_{k-l,1}+a_{k-l-1,l+3}+\ldots+a_{1,k+1}+a_{k+1,k+2}+
\ldots+a_{t,t+1}$$
if $l\leq k-2$ and
$$a=a_{1,1}+\ldots+a_{k-2,k-2}+a_{k+1,k+2}+
\ldots+a_{t,t+1}$$
if $l=k-1$.

Here $F_1,F_2$ and $G_1,G_2$ are two regular sequences with $F_1,F_2,G_1,G_2$ 
generic of degrees $a_{k-1,l},a_{k-1,l+1}+a_{k,l+2},a_{k,l+1},a_{k,l+2}$ passing through a common point, 
that is the common zero of the linear forms $L_1,L_2,L_3$ (see also Remark~\ref{defmod}). 
In particular, $$\alpha(\M_C)=a_{k-2,l-1}+\ldots+a_{k-l,1}+a_{k-l-1,l+3}+\ldots+a_{1,k+1}+a_{k+1,k+2}+
\ldots+a_{t,t+1}+1$$
if $l\leq k-2$ and
$$\alpha(\M_C)=a_{1,1}+\ldots+a_{k-2,k-2}+a_{k+1,k+2}+
\ldots+a_{t,t+1}+1$$
if $l=k-1$.
\end{rem}

\begin{rem}\label{smoothS}
If we do not require connectedness of $C$, we can perform the construction of Theorem~\ref{mmain} in such a
way that we have a surface
$S$ containing $C$ of degree\newline 
$a_{1,1}+\ldots+a_{t,t}+a$, for each $a>0$. $S$ can be taken
to be smooth on the complement of a zero-dimensional subset of $C$. Moreover, $S$ can be chosen in such a
way that its image in $I_X$ is a multiple of a minimal generator of minimal degree by a form of degree
$a>0$.
\end{rem}

\begin{proof}
In fact, for $t=2$, let $C=CI(F,G)\cup CI(H,J)$ be the disjoint union of
two generic, smooth, integral complete intersections. We have $deg\; F =a_{1,1}$, $deg\; G
=a_{1,2}+a_{2,3}$,
$deg\; H =a_{2,2}$, $deg\; J =a_{2,3}$. Then $C$ is smooth and contained in the surface of equation
$T=FH$.
$T$ has degree $a_{1,1}+a_{2,2}$, and its image in $I_X$ is a minimal generator. Let $S$ be the union of
$T$ with a generic surface $U$ of degree $a$. The singular locus of $T$ is $F\cap H$, so it is disjoint
from $C$. Let $Sing(S)$ denote the singular locus of $S$. $Sing(S)\cap C\subseteq U\cap C$, so it is a
zero-dimensional subset of $C$, by generality of $U$. The image of $S$ in $I_X$ is a multiple of the
minimal generator $T$ of minimal degree by the form $U$ of degree $a>0$.

Proceeding by induction on $t$, let $C=D\cup C_t$ be a basic double link of $D$ on a surface $S_1$ of
degree $a_{1,1}+\ldots+a_{t,t}$, with a general form of degree $a_{t,t+1}$. By the induction hypothesis
applied to $D$, we can choose $S_1$ smooth on the support of $D$, except possibly
for a zero-dimensional subset. By generality of our choice of the form of degree
$a_{t,t+1}$, we can also assume that the surface individuated by this form does not pass through any of the
singular points of $S_1$ contained in $D$. Let
$X,Y$ be the general plane sections of $C,D$ respectively. We can assume that the image of $S_1$ in $I_Y$
is a multiple of a minimal generator of minimal degree by a form of degree $a_{t,t}>0$. The image of
$S_1$ in $I_X$ is a minimal generator, by construction. Let
$S=S_1\cup U$,
$U$ a generic surface of degree
$a$. By generality of $U$, we can assume that $U$ does not pass through any of the points of $D\cap
C_t$ and that $U\cap C$ is zero-dimensional.
$Sing(S)=Sing(S_1)\cup (S_1\cap U)$, so $Sing(S)\cap C=(Sing(S_1)\cap D)\cup (Sing(S_1)\cap C_t)\cup
(S_1\cap U\cap C).$ $Sing(S_1)\cap D$ is zero-dimensional by assumption,
$Sing(S_1)\cap C_t$ is zero-dimensional, since $Sing(S_1)\cap C_t\cap D$ is empty by assumption.
$S_1\cap U\cap C$ is zero-dimensional, since $U\cap C$ is.
The image of $S$ in $I_X$ is a multiple of the minimal generator of
minimal degree image of $S_1$ by a form of degree $a>0$ (the image of $U$).
This is the proof, in the case $a_{t-1,1}\leq 2$. The proof in the other cases (see the Proof of
Theorem~\ref{mmain}) are analogous. 
\end{proof}

\begin{rem}\label{gens}
$I_C$ as constructed in Theorem~\ref{mmain} or in Remark~\ref{smoothS} is minimally generated in degree
less than or equal to
$$a_{k-1,l+1}+2a_{k,l+2}+a_{k-2,l-1}+\ldots+a_{k-l,1}+a_{k-l-1,l+3}+\ldots+a_{1,k+1}+a_{k+1,k+2}+
\ldots+a_{t,t+1}$$ $$=a_{1,2}+\ldots+a_{t,t+1}+a_{l-1,1}-a_{l-1,l}+a_{k,l+2}.$$
Notice that
$a_{1,2}+\ldots +a_{t,t+1}$ is the highest degree of a minimal generator of $I_X$.

One can easily show it proceeding by induction on $t$, and using the short exact sequence 
$$0\lra S(-s-t)\lra I_D(-t)\oplus S(-s)\lra I_C\lra 0$$
connecting the ideal of a scheme $D$ with
the ideal of its basic double link $C$ on a surface $S$ of degree $s$, with a form $F$ of degree $t$
(see Proposition 5.4.5 in \cite{M}). The case of a degree matrix of size $2\times 3$ is examined in
Example~\ref{2x3case}, and can be used as the basis of the induction.

As in Remark~\ref{components}, one can give a simple description of the deficiency module of the curves 
constructed in Remark~\ref{smoothS}.

In general, starting the construction from different submatrices of $M$ will yield curves that are not
projectively isomorphic. In fact, they are unions of complete intersections of different degrees.
The observations of Remark~\ref{smoothS} remain true, since they are independent
of which submatrix we start the construction from.
\end{rem}

\begin{rem}\label{nice}
If we start the construction of Theorem~\ref{mmain} from one of the curves constructed in 
Theorem~\ref{2x3nice}, we obtain a curve $C$ whose saturated ideal $I_C$ is generated in degree smaller
than or equal to $$a_{k-1,l+1}+a_{k,l+2}+1+a_{k-2,l-1}+\ldots+a_{k-l,1}+a_{k-l-1,l+3}+\ldots+a_{1,k+1}+
a_{k+1,k+2}+\ldots+a_{t,t+1}=$$ $$a_{1,2}+\ldots+a_{t,t+1}+a_{l-1,1}-a_{l-1,l}+1$$ if $l\leq k-2$, 
and in degree less than or equal to 
$$a_{k-1,k}+a_{k,k+1}+1+a_{k-2,k-2}+\ldots+a_{1,1}+a_{k+1,k+2}+
\ldots+a_{t,t+1}=$$ $$a_{1,2}+\ldots+a_{t,t+1}+a_{k,1}-a_{k,k-1}+1$$ if $l=k-1$.

$C$ is again a union of $t$ complete intersections. The same considerations as in Remark~\ref{components} 
about writing $C$ explicitly as a union of complete intersections hold. Using Remark~\ref{MCnice} and 
the Hartshorne-Schenzel Theorem, one can compute explicitly the initial and final degrees of the deficiency 
module of the curve, in terms of the entries of the degree matrix $M$ of its general plane section.
\end{rem}

\begin{rem}\label{singlocusC}
The space curve $C$ that we constructed in Remark~\ref{nice} is reduced, connected, non-degenerate, and 
non-aCM. We can take the complete intersections that constitute $C$ to be smooth,
so that $C$  has singularities only at the points of intersections of its irreducible components.
\end{rem}

We now find smooth surfaces that contain the curves constructed in Remark~\ref{nice}.
They will be used in the following constructions.

\begin{lemma}\label{smoothnice}
Let $C\subset\PP^3$ be a curve as constructed in Remark~\ref{nice}. Assume that the saturated
ideal of $C$ is minimally generated in degree smaller than or equal to $d$. Then there is a smooth
surface of degree $d$ containing $C$.
\end{lemma}

\begin{proof}
We proceed by induction on the size $t$ of the degree matrix of a general plane section of $C$. 
The case $t=2$ has been proved in Lemma~\ref{smnice}.

Consider the linear system $\Delta$ of surfaces of $\PP^3$ of degree $d$ containing $C$.
$C=C_1\cup C_2\cup\ldots\cup C_t$ is a reduced union of $t$ complete intersection curves. 
Let $Sing(C)=\cup_{i<j} C_i\cap C_j$ be the singular locus of $C$ (see
Remark~\ref{singlocusC} about what the singular locus of $C$ looks like). 
The general element of $\Delta$ is basepoint-free outside of $C$, hence smooth
outside of $C$ by Bertini's Theorem. Consider now a point $P\in C$. We want to show that the 
general element of $\Delta$ is smooth at $P$. By Corollary 2.10 in \cite{GV},
it is enough to exhibit two elements of $\Delta$ meeting transversally at $P$. Since $C$
is smooth outside of $Sing(C)$, for each point $P\not\in Sing(C)$ we have two minimal generators of
$I_C$, call them $F$ and $G$, meeting transversally at $P$. The degree of each of them is at most $d$. 
Add generic planes as needed, to obtain surfaces of degree $d$ that meet transversally at $P$. 
In order to complete the proof, we need to check that the points of $Sing(C)$ are not fixed singular
points for $\Delta$. So it is enough to find a surface for each $P\in Sing(C)$ that contains $C$ and is
non-singular at $P$. Each singular point of $C$ is the intersection of two irreducible components 
of the curve, $P\in C_i\cap C_j$ for some $1\leq i<j\leq t$. We can assume, 
by generality of our choices, that $i,j$ are determined by $P$, i.e. we can assume that there are exactly 
two irreducible components of $C$ meeting at $P$. 
Without loss of generality, we can then assume that $j=t$ and that 
$C=D\cup C_t$, where $P\not\in Sing(D)$. As seen in Remark~\ref{smoothS}, we can perform the 
basic double link in such a way that the surface $S_1$ of degree $a_{1,1}+\ldots+a_{t,t}$ that we perform 
the link on is smooth on $D$ outside of a zero-dimensional subscheme. Moreover, we can assume that the 
singular locus of $S_1$ does not contain any of the points of $D\cap C_t$. 
In particular, $S_1$ is smooth at $P$ and contains $C$. Notice that 
$d\geq a_{1,1}+\ldots+a_{t,t}=\alpha(I_C)$. Add to $S_1$ a generic surface of degree 
$d-a_{1,1}-\ldots-a_{t,t}$ to obtain a surface containing $C$ and smooth at $P$.
\end{proof}

We can also ask the question whether it is possible to give a sufficient condition for the 
Cohen-Macaulayness of $C\subset\PP^3$, in terms of the entries of the $h$-vector of its general 
plane section $X\subset\PP^2$. It is easy to see that we cannot, as the following proposition shows.

\begin{prop}\label{vect}
Let $h(z)=1+h_1z+\ldots+h_sz^s$, $h_s\neq 0$ be the $h$-vector of some zero-dimensional scheme in $\PP^2$.
Then there exists a non-aCM, reduced curve $C\subset\PP^3$, whose general plane section
$X\subset\PP^2$ has $h$-vector $h(z)$. Moreover, the curve $C$ can be taken to be connected, 
unless $h(z)=1+2z$.
\end{prop}

\begin{proof}
To any $h$-vector $h(z)$, we can uniquely associate a degree matrix $M$ with no entries equal to $0$, such
that if $X\subset\PP^2$ is a zero-dimensional scheme with degree matrix $M$, then the $h$-vector of $X$ is $h(z)$. 
If $M$ has one entry less than or equal to $2$ and is not a $2\times 3$ matrix with all its entries equal
to $1$, by Theorem~\ref{mmain} we can find a non-aCM, reduced, connected curve $C\subset\PP^3$, whose 
general plane section $X\subset\PP^2$ has degree matrix $M$, hence $h$-vector $h(z)$. 
If $M$ is the degree matrix of size $2\times 3$ with all entries equal to $1$, i.e. if the $h$-vector is
$h(z)=1+2z$, let $C$ be the disjoint union of a reduced plane conic and a line.  
If $M=(a_{i,j})_{i=1,\ldots,t;\; j=1,\ldots, t+1}$ has $a_{t,1}\geq 3$, let
$N=(b_{i,j})_{i=1,\ldots,t+1;\; j=1,\ldots,t+2}$ be the degree matrix with entries $b_{i,j}=a_{i,j-1}$ for
$i=1,\ldots,t$, $j=2,\ldots,t+2$, $b_{t+1,1}=0$, $b_{t+1,2}=2$. $N$ is determined by these entries, under
the assumption that it is homogeneous. $b_{i,j}>0$ for $(i,j)\neq (t+1,1)$, so $N$ is a degree matrix.
Moreover, the $h$-vector of a zero-dimensional scheme that has degree matrix $N$ is again $h(z)$. Then, by 
Theorem~\ref{mmain}, there exists a non-aCM, reduced, connected curve $C\subset\PP^3$, whose 
general plane section $X\subset\PP^2$ has degree matrix $M$, hence $h$-vector $h(z)$.  
\end{proof}

Let us now look at the general case of a curve $C\subset\PP^{n+1}$, whose general hyperplane section 
is the zero-dimensional scheme $X\subset\PP^n$. With the notation of Definition~\ref{liftingmatrix}, let $M=(a_{ij})$
be the lifting matrix of $X$.

We have a sufficient condition for the Cohen-Macaulayness of $C$, analogous to the case $n=2$,
that follows again from Theorem~\ref{socle}.

\begin{cor}\label{bigger3P^n}
Let $C\subset\PP^{n+1}$ be a curve, whose general hyperplane section $X\subset\PP^n$
has lifting matrix $M=(a_{ij})$. If $a_{t,1}\geq n+1$, then $C$ is arithmetically Cohen-Macaulay.
\end{cor}

\begin{proof}
With the notation of Theorem~\ref{socle}, if $C$ is not aCM, we have $$b\geq m_t-n\geq d_1+1.$$
Then all the minimal generators of $I_X$ lift to $I_C$, so $C$ is aCM.  
\end{proof}

\section{What can be said about integral curves?}

Throughout this section, we will concentrate on integral (reduced and
irreducible), locally Cohen-Macaulay, equidimensional, non-degenerate curves $C\subset\PP^3$. We want
to investigate whether, under the extra assumption of integrality on the curve, we can find a condition on
the degree matrix of $X$, that is weaker than in Corollary~\ref{bigger3}, and still forces $C$ to be
arithmetically Cohen-Macaulay.

First of all, we need a characterization of the matrices of integers that can occur as the degree matrix
of  a zero-dimensional scheme in $\PP^2$ that is the (general) plane section of an integral, aCM space curve. We will
call such a matrix $M$ an {\em integral degree matrix}. 
Homogeneous matrices of integers that can occur as integral degree matrices have been characterized by
J. Herzog, N.V. Trung and G. Valla in \cite{HTV}. In our language, they prove the following result.

\begin{thm}\label{htv}
Let $M=(a_{i,j})$ be a matrix of integers of size $t\times (t+1)$. Then $M$ is an integral degree matrix
if and only if it is homogeneous and $a_{h+1,h}>0$ for $h=1,\ldots,t-1$.
\end{thm}

We will start our analysis by looking at an example.

\begin{ex}\label{exirred}
Consider the following degree matrix
$$M=\left(\begin{array}{ccc} 1 & 3 & 3 \\ 1 & 3 & 3 \end{array}\right)$$
corresponding to some zero-dimensional scheme $X$ of degree $deg(X)=15$. $I_X$ has minimal free resolution
$$0\lra R(-7)^2 \lra R(-6)\oplus R(-4)^2 \lra I_X \lra 0.$$
The construction of Theorem~\ref{main} yields an example of a reduced, connected space curve, that is
non-aCM and such that its general plane section has degree matrix $M$. $X$ is in fact the general
plane section of $C=CI(1,6)\cup CI(3,3)$, where the two complete intersections are generic, through a
common point.

Assume now that we have a reduced, irreducible curve $C\subset \PP^3$ whose general plane section $X$
has degree matrix $M$. By Theorem~\ref{socle}, the
minimal degree of an element of $I_X$ that is not the image of some element of $I_C$ under the quotient
map is $b\geq 7-2=5$.
Then the two minimal generators of $I_X$ of degree 4 are the images of two minimal generators $F,G$
of $I_C$.
Moreover, since $C$ is integral, $F,G$ are both irreducible forms. Hence they form a regular
sequence in 
$S=k[x_0,x_1,x_2,x_3]$. Let $E$ be the curve whose saturated ideal is $I_E=(F,G)\subset S$.
$E$ is a complete intersection and it contains $C$, hence $C$ is linked via $E$ to a curve $D$.
$D$ has degree $deg(D)=deg(E)-deg(C)=1$ (see \cite{M}, Corollary 5.2.13), so it is a line. In
particular, $D$ is aCM. Since the property of being aCM is an invariant of the CI-linkage class
of a scheme (see \cite{M}, Theorem 5.3.1), $C$ has to be aCM as well.
\end{ex}

The example inspires the following observations.

\begin{lemma}\label{lift}
Let $C\subset \PP^3$ be a curve, whose general plane section $X$ has degree matrix $M=(a_{i,j})$ of size 
$t\times (t+1)$.
Assume that $a_{t,j}\geq 3$. Then the $t+2-j$ minimal generators of lowest degrees of $I_X$ are images 
of the $t+2-j$ minimal generator of lowest degrees of $I_C$.
\end{lemma}

\begin{proof}
It directly follows from Theorem~\ref{socle}. Let $d_j,\ldots,d_{t+1}$ be the degrees of the $t+2-j$ 
minimal generators of lowest degrees of $I_X$, $d_{t+1}\leq\ldots\leq d_j$ (here we follow the 
notation of Theorem~\ref{socle}; notice that some of the degrees could be repeated). 
The lowest shift in the last free module of the minimal free resolution of $I_Z$ is 
$d_{t+1}+a_{t,t+1}=d_j+a_{t,j}$. 
 
If $a_{t,j}\geq 3$, by Theorem~\ref{socle} it follows that the minimum degree of a polynomial in $I_X$ 
that is not the image of an element of $I_C$ under the standard projection map is 
$b\geq d_j+a_{t,j}-2>d_j$. Therefore the $t+2-j$ minimal generators of lowest degrees of $I_X$ 
are images of minimal generators of $I_C$.
\end{proof}

We can now state the first condition that forces an integral curve $C$ to be aCM. The condition is given 
in terms of the entries of the degree matrix of its general plane section. The proof of the Proposition
is a generalization of the argument of Example~\ref{exirred}.

\begin{prop}\label{aCM2x3}
Let $C\subset \PP^3$ be a reduced, irreducible curve, whose general plane section $X$ has degree
matrix
$M=(a_{i,j})$ of size $2\times 3$. Assume that $a_{2,2}\geq 3$ and that $a_{1,1},a_{2,1}\not= 2$.
Then $C$ is aCM.
\end{prop}

\begin{proof}
Under these assumptions, it follows from Lemma~\ref{lift} that the two generators of minimal degrees
of $I_X$
lift to two minimal generators of $I_C$, call them $F,G$. Following the strategy of Example~\ref{exirred},
we  notice that $F,G$ form a regular sequence in $S$.
Let $E$ be the complete intersection corresponding to $I_E=(F,G)\subset S$. Let $D$ be the curve
residual to $C$ 
in the link. Taking general plane sections, the link is preserved and we obtain that the general
plane 
section $Y$ of 
$D$ has degree matrix $(a_{1,1},a_{2,1})$. By the already-mentioned result of Strano (Theorem 6,
\cite{S}), $D$ has to be aCM and, 
since the property of being aCM is an invariant of the CI-linkage class of a scheme 
(see \cite{M}, Theorem 5.3.1), $C$ has to be aCM as well.
\end{proof}

In what follows, we will make extensive use of Bertini's Theorem. For our convenience, we recall it
here in the form we will need it. See \cite{H}, Ch. III, Corollary 10.9 and the following remark for a
proof.

\begin{thm}(Bertini)\label{bert}
Let $S$ be a (smooth) integral projective scheme of dimension at least $2$, over an
algebraically closed field of characteristic $0$. Let $\delta$ be a basepoint-free linear
system on
$S$. Then a generic element of $\delta$ is a (smooth) integral subscheme of $S$.
\end{thm}

Using Bertini's Theorem, we can find another family of degree matrices $M$ such that every integral space 
curve $C$ whose general plane section $X$ has degree matrix $M$ is arithmetically Cohen-Macaulay.

\begin{prop}\label{aCM-II}
Let $C\subset \PP^3$ be a reduced, irreducible curve, whose general plane section $X$ has 
degree matrix
$M=(a_{i,j})$ of size $2\times 3$. Assume that $a_{1,1},a_{2,3}\geq 3$ and that 
$a_{2,1},a_{2,2}\not= 2$. Then $C$ is aCM.
\end{prop}

\begin{proof}
Since $a_{2,3}\geq 3$, we can conclude by Lemma~\ref{lift} that the generator of minimal 
degree of $I_X$ lifts to a minimal generator of $I_C$. Then, $I_X$ and $I_C$ have the
same 
initial degree $\alpha=a_{1,1}+a_{2,2}$. Let $T$ be a surface of degree $\alpha$ 
containing $C$. $T$ is integral since $C$ is integral. 

Consider the linear system $\Sigma_d$ on $T$ of the 
curves cut out on $T$, outside of $C$, by the surfaces of degree $d$ containing $C$.
For $d\gg 0$, in particular for $d$ bigger or equal to the highest degree of a minimal 
generator of $I_C$, the linear system $\Sigma_d$ is basepoint-free. By Bertini's 
Theorem, its general element is an integral curve, call it $D$. $D$ is CI-linked to $C$
by construction, so its general plane section is CI-linked to the general plane 
section of $C$ via a $CI(\alpha,d)$.
Let $Y$ be the general plane section of $D$. The degree matrix of the general 
plane section $X$ of $C$ is $M$, hence (see \cite{M}, Proposition 5.2.10) a minimal free resolution for
$I_Y$ is 
$$\begin{array}{ccccc}
 & R(-d-a_{1,1}+a_{1,3}) & & R(-a_{1,1}-a_{2,2})\oplus R(-d+a_{1,3}) & \\
0\lra & \oplus & \lra & \oplus & \lra I_Y\lra 0 \\
 & R(-d-a_{2,2}+a_{2,3}) & & R(-d+a_{2,3}) & 
\end{array}$$
since the form of degree $\alpha$ is a minimal generator of $I_X$, while the form of degree $d$ is not. 
Then the degree matrix of $I_Y$ is
$$N=\left( \begin{array}{ccc}
a_{2,2} & a_{1,2} & d-a_{1,1}-a_{2,3} \\
a_{2,1} & a_{1,1} & d-a_{2,2}-a_{1,3}
\end{array} \right).$$
Since $d\gg 0$, we can assume that $d-a_{2,2}-a_{1,3}\geq a_{1,1}$ (notice that this also guarantees
minimality of the resolution of $I_X$ above). By hypothesis we have
$a_{1,1}\geq 3$ and $a_{2,1},a_{2,2}\not= 2$, so we can apply Proposition~\ref{aCM2x3} to
conclude that $D$ is aCM. Then $C$ is aCM as well.
\end{proof}

\begin{rem}
In our situation, assuming $a_{2,1}\neq 2$ is equivalent to $a_{2,1}=1$.
\end{rem}

From Proposition~\ref{aCM2x3} and Proposition~\ref{aCM-II}, we can deduce some conditions on the 
$h$-vector of $X$ that force $C$ to be arithmetically Cohen-Macaulay.

For what follows we need to derive a formula for the $h$-vector of a zero-dimensional scheme 
$X\subset\PP^2$ in terms of the entries of the degree matrix of $X$. 

\begin{lemma}\label{mat-hvect}
Let $X\subset\PP^n$ be an arithmetically Cohen-Macaulay scheme of codimension $2$, and let
$M=(a_{i,j})_{i=1,\ldots,t;\; j=1,\ldots,t+1}$ be its degree matrix. Then the $h$-vector of $X$ is
$$h(z)=\sum_{i=1}^t z^{a_{1,1}+\ldots
+a_{i-1,i-1}}(1+z+\ldots+z^{a_{i,i}-1})(1+z+\ldots+z^{a_{i+1,i+1}+\ldots+a_{t,t}+a_{i,t+1}}).$$
\end{lemma}

\begin{proof}
The minimal free resolution of $X$ is 
$$0\lra \FF_2\lra \FF_1 \lra I_X\lra 0$$ 
where 
$$\FF_2=\oplus_{i=1}^t R(-a_{1,1}-\ldots-a_{t,t}-a_{i,t+1}),$$
$$\FF_1=\oplus_{j=1}^t R(-a_{1,1}-\ldots-a_{t,t}+a_{j,j}-a_{j,t+1})\oplus R(-a_{1,1}-\ldots-a_{t,t}),$$ 
$I_X\subset
R=k[x_0,\ldots,x_n]$. Then the $h$-vector of
$X$ is $$h(z)=\frac{1-\sum_{i=1}^t
z^{a_{1,1}+\ldots+a_{t,t}-a_{i,i}+a_{i,t+1}}-z^{a_{1,1}+\ldots+a_{t,t}}+\sum_{i=1}^t
z^{a_{1,1}+\ldots+a_{t,t}+a_{i,t+1}}}{(1-z)^2}.$$
Computing, we get $$1-z^{a_{1,1}+\ldots+a_{t,t}}+\sum_{i=1}^t 
(z^{a_{1,1}+\ldots+a_{t,t}+a_{i,t+1}}-z^{a_{1,1}+\ldots+a_{t,t}-a_{i,i}+a_{i,t+1}})=$$
$$=(1-z)[(1+z+\ldots+ z^{a_{1,1}+\ldots+a_{t,t}-1})-\sum_{i=1}^t
z^{a_{1,1}+\ldots+a_{t,t}-a_{i,i}+a_{i,t+1}}(1+z+\ldots+z^{a_{i,i}-1})]=$$ $$=(1-z)^2 \sum_{i=1}^t
z^{a_{1,1}+\ldots+a_{i-1,i-1}}
(1+z+\ldots+z^{a_{i,i}-1})(1+z+\ldots+z^{a_{i+1,i+1}+\ldots+a_{t,t}+a_{i,t+1}}).$$
\end{proof}

\begin{rem}
The degree matrix of a scheme $X$ as in Lemma~\ref{mat-hvect} determines the $h$-vector of
$X$, while the $h$-vector of $X$ determines the degree matrix only under the hypothesis that all the
entries of the degree matrix of $X$ are positive.
\end{rem}

\begin{rem}
From Lemma~\ref{mat-hvect}, we see that the $h$-vector of $X$ can be formally written as a sum of some
shifts of the $h$-vectors $h_i(z)$ of $t$ complete intersections of type
$(a_{i,i},a_{i+1,i+1}+\ldots+a_{t,t}+a_{i,t+1}+1)$, $i=1,\ldots,t$. 
The $h$-vector $h_i(z)$ has increasing coefficients in degrees $1,\ldots,a_{i,i}-1$, they are
constant until degree $a_{i+1,i+1}+\ldots+a_{t,t}+a_{i,t+1}$, and then they are decreasing. 
Looking at $k_i(z)=z^{a_{1,1}+\ldots+a_{i-1,i-1}}h_i(z)$, we have that the coefficients start decreasing
in degree $f_i=a_{1,1}+\ldots+a_{i-1,i-1}+a_{i+1,i+1}+\ldots+a_{t,t}+a_{i,t+1}$ and the last nonzero
coefficient appears in degree $e_i=f_i+a_{i,i}-1$.

Under the assumption that the degree matrix $M$ is integral, we have $a_{i+1,i}>0$ for all $i$, that
gives $e_{i+1}-f_i=f_{i+1}+a_{i+1,i+1}-1-f_i=a_{i+1,i}-1\geq 0$, so each $k_{i+1}(z)$ does not end
on the flat part of $k_i(z)$. 

This shows that the $h$-vector of $X$ is of decreasing type. Moreover, $h_j-h_{j+1}\geq 2$ for all
$f_i\leq j\leq e_{i+1}$ for some $i$, and only for those $j$'s.  
\end{rem}

We are now ready to derive some sufficient conditions for $C$ integral to be aCM, in terms of the
$h$-vector of its general plane section.

\begin{cor}\label{h-vect}
Let $C\subset \PP^3$ be a reduced, irreducible curve, whose general plane section $X$ has
$h$-vector $h=1+h_1 z+\ldots+h_s z^s$, $h_s\neq 0$.
Let $$u=\mmax\{ i\; |\; h_i=i+1\}, \;\; v=\mmax\{ i\; |\; h_i=u+1\},$$ $$w=\mmin\{ i\; |\; v\leq i\leq
s-1, \; h_i-h_{i+1}\neq 1\}.$$ If \/ $\{ i\; |\; v\leq i\leq s-1, \;
h_i-h_{i+1}\neq 1\}=\{ w \}$ and either
$$s=u+v-1,\;\; u+v-w\neq 2, \;\;\mbox{and}\;\; w-v\geq 2$$ or $$s=u+v-1,\;\; v\geq 6,\;\; w-u\geq
3,\;\;\mbox{and}\;\; w\neq v+1$$ then $C$ is arithmetically Cohen-Macaulay.
\end{cor}

\begin{proof}
Let $M=(a_{i,j})_{i=1,\ldots,t;\; j=1,\ldots,t+1}$ be the degree matrix of $X$. By
Lemma~\ref{mat-hvect}, the $h$-vector of
$X$ is $$h(z)=\sum_{i=1}^t z^{a_{1,1}+\ldots
+a_{i-1,i-1}}(1+z+\ldots+z^{a_{i,i}-1})(1+z+\ldots+z^{a_{i+1,i+1}+\ldots+a_{t,t}+a_{i,t+1}}).$$
Let
\begin{equation}\label{hvect}
h_i(z)=(1+z+\ldots+z^{a_{i,i}-1})(1+z+\ldots+z^{a_{i+1,i+1}+\ldots+a_{t,t}+a_{i,t+1}}).
\end{equation}
By Lemma~\ref{mat-hvect}, we can think of $h(z)$ as the sum of $t$ $h$-vectors of complete intersections
of types $(a_{i,i},a_{i+1,i+1}+\ldots+a_{t,t}+a_{i,t+1}+1)$, for $i=1,\ldots,t$. 
By assumption, $\{ i\; |\; v\leq i\leq s-1, \; h_i-h_{i+1}\neq 1\}=\{ w \}$, so $h_i-h_{i+1}=1$ for
$v\leq i\leq s-1$, $i\neq w$, so the $h$-vector has only one jump of more than $1$, once it starts
decreasing. Therefore, it has to be the sum of only two $h$-vectors $h_i$, that is $t=2$.
The degree matrix of $X$ has then size $2\times 3$. $X$ is the general
plane section of an integral curve $C$ (so it has UPP, see \cite{H} about the general plane section of an
integral curve and its $h$-vector). Then $M$ is integral, so in particular $a_{2,1}>0$. All the
entries of $M$ are positive, so the $h$-vector of $X$ determines the degree matrix.
From equation (\ref{hvect}), we can compute 
$$u=a_{1,1}+a_{2,2}-1,\;\;\;\;
v=a_{1,1}+a_{2,3},$$
\begin{equation}\label{equ} w={a_{1,1}+a_{2,2}+a_{2,3}-a_{2,1}} \;\;\;\mbox{and}\;\;\;
s=2a_{1,1}+a_{2,2}+a_{2,3}-2.\end{equation}  
The assumption that $\{ i\; |\; v\leq i\leq s-1, \; h_i-h_{i+1}\neq 1\}=\{ w \}$ forces $a_{2,1}=1$: in
fact, $h_i=h_{i+1}-2$ for
$w=a_{1,1}+a_{2,2}+a_{2,3}-a_{2,1}\leq i\leq a_{1,1}+a_{2,2}+a_{2,3}-1$. Solving the equations
(\ref{equ}) gives $$s=u+v-1$$ and $$a_{1,1}=u+v-w,\;\;\; a_{2,2}=w-v+1, \;\;\; a_{2,3}=w-u,$$ so the
degree matrix of
$X$ has the following form, in terms of $u,v,w$
$$M=\left(\begin{array}{ccc}
u+v-w & u & v-1 \\
1 & w-v+1 & w-u
\end{array}\right).$$
By Proposition~\ref{aCM2x3}, $C$ is aCM if $u+v-w\neq 2$ and $w-v\geq 2$.
By Proposition~\ref{aCM-II}, $C$ is aCM if $u+v-w\geq 3$, $w-u\geq 3$ and $w-v+1\neq 2$, or equivalently
if $w-u\geq 3$, $v\geq 6$ and $w\neq v+1$.
\end{proof}

For any degree matrix that has at least one entry smaller than or equal to $2$ and that
does not fall in one of the two classes of examples of
Proposition~\ref{aCM2x3} and Proposition~\ref{aCM-II}, 
we can produce an integral, smooth curve that is non-aCM, and whose general plane section
has degree matrix $M$.

The following lemmas will be needed for the construction of a smooth, integral curve whose general 
plane section has a prescribed degree matrix.

\begin{lemma}\label{smooth}
Let $C\subset\PP^3$ be a smooth space curve, whose ideal is minimally generated in degree
smaller than or equal to $d$. Then there is an integral, smooth surface of degree $d$
containing $C$.
\end{lemma}

\begin{proof}
Consider the linear system $\Delta$ of surfaces of $\PP^3$ of degree $d$, containing $C$.
The general element of $\Delta$ is basepoint-free outside of $C$, hence smooth outside of
$C$ by Bertini's Theorem. Consider now a point $P\in C$. By Corollary 2.10 in \cite{GV},
it is enough to exhibit two elements of $\Delta$ meeting transversally at $P$. Since $C$
is smooth, for each point we have two minimal generators of $I_C$, call them $F$ and $G$,
meeting transversally at $P$. The degree of each of them is at most $d$. Add generic
planes as needed, to obtain surfaces of degree $d$ that meet transversally at $P$. 

A general surface of degree $d$ containing $C$ will be integral as well. In fact, if it
were not, all the minimal generators of $I_C$ would have to share a factor, but that is
not possible since $I_C$ has height 2.  
\end{proof}

\begin{lemma}\label{deform}
Let $C$ be a curve lying on a smooth surface $S\subset\PP^3$. Let $M$ be the degree matrix of the
general plane section of $C$, and assume that all the entries of $M$ are different from zero. Let $D$ be a
general element of the linear system $|C|$. Then the degree matrix of the general plane section of $D$ is
$M$.   
\end{lemma}

\begin{proof}
Let $X,Y\subset\PP^2$ be the general plane sections of $C,D$ respectively. We have $\I_X\cong \I_Y$
as $\mathcal{K}$-modules, where $\mathcal{K}$ the sheaf of total quotients of
the structure sheaf of $S\cap H$ for a general plane $H=\PP^2$. Then $X$ and $Y$ have the same
Hilbert function (their graded parts have the same dimension as $H^0(\mathcal{K})$-vector spaces,
hence as $k$-vector spaces). Notice that $H^0(\mathcal{K})$ is a field, since $S\cap H$ is an
integral curve by Bertini's Theorem.
Look now at the linear system $|X|$ of divisors on $S\cap H$ that are linearly equivalent to $X$.
$Y$ is the general element of $|X|$ and since the degree matrix of $X$ has no zero entries, neither
does the degree matrix of $Y$, by upper-semicontinuity. Then $X,Y\subset\PP^2$ both have
degree matrices with non-zero entries and the same Hilbert series, so they have the same degree matrix.
\end{proof}

We are now ready to construct a smooth, integral, non-aCM curve, whose general plane section has
a prescribed degree matrix. We can perform the construction for each integral matrix
such that at least one of the entries is smaller than $3$ and that does not fall in the classes of
examples covered by Proposition~\ref{aCM2x3} and Proposition~\ref{aCM-II}. We exclude from our
analysis the degree matrix of size $2\times 3$ with all entries equal to $1$ (see Example~\ref{mat1}).

We will start with an analysis of the degree matrices of size $2\times 3$.

\begin{thm}\label{smooth2x3}
Let $M=(a_{i,j})_{i=1,2; \; j=1,2,3}$ be a degree matrix with positive entries, such that
$a_{2,1}\leq 2$. Suppose that the entries of $M$ are not all equal to $1$, and that they do not satisfy
the hypothesis of either Proposition~\ref{aCM2x3}, or Proposition~\ref{aCM-II}. 
Then there exists an integral, smooth, non-aCM curve in
$\PP^3$ whose general plane section has degree 
matrix $M$.
\end{thm}

\begin{proof}
We'll be performing different constructions, depending on the entries of the
matrix $M$. Notice that, for $2\times 3$ matrices, being integral is equivalent to having positive
entries, since $a_{2,1}>0$. 
 
{\bf Case 1.} Assume $a_{2,1}=2$.

In this case $$M=\left( \begin{array}{ccc} 
a_{1,1} & a_{1,1}+a_{2,2}-2 & a_{1,1}+a_{2,3}-2 \\
2 & a_{2,2} & a_{2,3} 
\end{array} \right).$$
Let $D$ be a general rational smooth curve of degree $2a_{1,1}$, lying on a smooth quadric
surface. Taking $D$ as in Remark~\ref{rat2}, we can assume that the saturated ideal $I_D$ is
generated in degree less than or equal to $a_{1,1}+2$. The general plane section of $D$ is a Complete
Intersection of type $(2,a_{1,1})$. 
Let $F$ be the equation of a smooth, integral surface of degree $a_{1,1}+a_{2,2}$ containing $D$. Such an
$F$ exists, by Lemma~\ref{smooth}. Consider the linear system of curves cut out on $F$ outside of $D$ by
surfaces of degree $a_{1,1}+a_{2,3}$ containing $D$. The linear system is base-point free, since
$a_{1,1}+a_{2,3}\geq a_{1,1}+2$. So by Bertini's Theorem, the general element $C$ is a smooth, integral
curve. By construction, $C$ is linked to $D$ via a $CI(a_{1,1}+a_{2,2},a_{1,1}+a_{2,3})$.
Then, by Proposition 5.2.10 in \cite{M}, the general plane section $X$ of $C$ has a free resolution
$$\begin{array}{ccccc} & R(-a_{1,1}-a_{2,2}-a_{2,3}) & & R(-a_{1,1}-a_{2,2}-a_{2,3}+2) & \\
0\ra & \oplus & \ra & \oplus & \ra I_X\ra 0. \\
 & R(-2a_{1,1}-a_{2,2}-a_{2,3}+2) & & R(-a_{1,1}-a_{2,2})\oplus R(-a_{1,1}-a_{2,3}) & \end{array}$$
Then the general plane section $X$ of $C$ has degree matrix $M$. No cancellation can occur, 
since all the entries of $M$ are positive. Hence the free resolution is minimal and $X$ has degree
matrix $M$.

{\bf Case 2.}  Assume $a_{2,1}=1$ and $a_{2,2}=2$.
 
Let $D$ be two skew lines, 
$$M=\left( \begin{array}{ccc} 
a_{1,1} & a_{1,1}+1 & a_{1,1}+a_{2,3}-1 \\
1 & 2 & a_{2,3} 
\end{array} \right).$$

Perform a basic double link using generic polynomials $F\in I_D$ and
$G\in S=k[x_0,x_1,x_2,x_3]$, of  degrees $a_{1,1}+2$ and $a_{1,1}+a_{2,3}-1$ respectively. 
We obtain a curve $C=D\cup (F\cap G)$, whose general plane section has degree matrix
$M$. In fact, we have the exact sequence (see \cite{M}, Theorem 3.2.3 and Remark 3.2.4 b)
$$0\lra R(-2a_{1,1}-a_{2,3}-1)\lra I_{D\cap H}(-a_{1,1}-a_{2,3}+1)\oplus
R(-a_{1,1}-2) \lra I_{C\cap H}\lra 0$$ 
for $H$ a general plane in $\PP^2$, $R=k[x_0,x_1,x_2]$.

The surface defined by $F$ is smooth and integral by generality, and the linear system
$|C|$ of curves on $F$ that are linearly equivalent to $C$ is  basepoint-free. 
In fact, the linear system $|D|$ on $F$ is itself basepoint-free: let $P$ be a point
of $D$ and let $U$ be the equation of a generic surface of degree $d$ containing $D$.
For $d\gg 0$, $U$ will be the equation of a smooth surface, containing $D$ and meeting
$F$ transversally. Let $U\cap F=D\cup D'$.
By generality, we can assume that $P\not\in D'$. Let
$T$ be the equation of a generic surface of the same degree
$d$, containing the curve $D'$.
$F\cap T=D'\cup E'$. By generality assumption the surface $T$, hence the curve $E'$, does
not pass through $P$ and the divisor
$D-(U\cap F)+(T\cap F)=D-D-D'+D'+E'=E'$ is linearly equivalent to $D$.
Hence, $|D|$ is basepoint-free.
By Bertini's Theorem, $|C|$ contains a 
smooth, integral, non-aCM curve, whose general plane section has degree matrix $M$ by Lemma~\ref{deform}.

{\bf Case 3.} Assume $a_{2,1}=1$ and $a_{1,1}=2$.

In this case, the degree matrix $M$ is
$$\left( \begin{array}{ccc}
2 & a_{1,2} & a_{1,3} \\
1 & a_{1,2}-1 & a_{1,3}-1
\end{array} \right).$$
Let $D$ be two skew lines. Its general plane section consists of two distinct points, hence it has 
degree matrix $(1,2)$. Let $U$ be a smooth surface of degree $a_{1,2}+1\geq 3$ containing
$D$. Let $C$ be the general element of the 
linear system cut out on $U$, outside of $D$, by the surfaces of degree $a_{1,3}+1\geq 3$.
The linear system is basepoint-free outside of $D$, since the ideal $I_D$ is generated 
entirely in degree $2$. The general element of the linear system links $D$ to the curve $C$,
that is smooth and integral by Bertini's Theorem. Moreover, $C$ is not arithmetically
Cohen-Macaulay, since $D$ is not.

The general plane sections $X,Y$ of $C,D$ are
CI-linked  via a complete intersection of type $(a_{1,2}+1,a_{1,3}+1)$. By
Proposition 
5.2.10 in \cite{M} we have the following free resolution for $X$:
$$\begin{array}{rcccl}
& R(-a_{1,2}-a_{1,3}) & & R(-a_{1,2}-1)\oplus R(-a_{1,3}-1) & \\
0\lra & \oplus  & \lra & \oplus & \lra I_X\lra 0. \\
& R(-a_{1,2}-a_{1,3}-1) & & R(-a_{1,2}-a_{1,3}+1) & 
\end{array}$$
So the degree matrix of the general plane section $X$ of $C$ is $M$. No cancellation can 
occur in the free resolution of $X$, since none of the entries of $M$ is zero.

{\bf Case 4.} Assume $a_{2,1}=1$ and $a_{1,1}=1$.

By Proposition~\ref{aCM2x3} $a_{2,2}\leq 2$, hence we can assume $a_{2,2}=1$ (the situation
when $a_{2,2}=2$ is treated in Case 2). The degree matrix is then of the form
$$M=\left( \begin{array}{ccc}
1 & 1 & a \\
1 & 1 & a
\end{array} \right)$$
for some $a\geq 2$. For $a=1$, assuming $C$ integral or even $C$ connected, already forces $C$ to be aCM
(see Example~\ref{mat1}). If $a=2$, we can let $C$ be a general smooth rational curve of degree
$5$. Its general plane section consists of $5$ generic points in $\PP^2$, as we showed in
Example~\ref{genpoints}. Hence it has degree matrix $M$.

For any $a\geq 2$, let $D$ consist of $2a-1$ skew lines on a
smooth quadric surface $Q$. The general plane section $Y$ of $D$ has degree matrix 
$$N=\left( \begin{array}{ccc}
1 & 1 & a-1 \\
1 & 1 & a-1
\end{array} \right)$$
and $I_D$ is generated in degrees $2,a$. Let $E$ be the complete intersection whose
defining ideal is
$I_E=(Q,F)$. Here $F$ is the equation of a generic surface of degree $2a$ containing
$D$. Let $F$ vary among all the surfaces of degree $2a$
containing $D$. Consider the linear system of curves that are residual to $D$ in
the complete intersection $E$. Bertini's Theorem applies, since the linear system is base-point free. 
Then the residual $C$ to $D$ in $E$, for a generic $F$, is smooth,
integral and non-aCM.

Applying Proposition 5.2.10 in \cite{M} to the general sections $X,Y$ of $C,D$, we
get that the minimal free resolution of $X$ is
$$0\lra R(-a-2)^2\lra R(-a-1)^2\oplus R(-2)\lra I_X\lra 0,$$
hence $C$ is smooth, integral, non-aCM and its general plane section has degree matrix
$M$.

{\bf Case 5.} Assume $a_{2,1}=1$ and $a_{1,1}\geq 3$.

We can assume $a_{2,2}=1$, since
the case $a_{2,2}=2$ has been treated in Case 2. Moreover, by Proposition~\ref{aCM-II} we have 
$a_{2,3}\leq 2$. The proof in the case $a_{2,3}=2$ is analogous to the proof of Case 2, starting with
$D$ equal to two skew lines and performing a basic double link using forms $F\in I_C$ and $G\in S$, of
degrees $a_{1,1}+1,a_{1,2}$ respectively.

Suppose then that $a_{2,3}=1$. Let $$N=\left(
\begin{array}{ccc} 
a & a & a \\
1 & 1 & 1
\end{array} \right)$$ and let $D$ consist of $2a+1$ skew lines on a smooth quadric
surface. The ideal $I_D$ is generated in degrees $2,a+1$, and the degree matrix of the
general plane section $Y$ of $D$ is
$$N=\left(
\begin{array}{ccc} 
1 & 1 & a \\
1 & 1 & a
\end{array} \right).$$
Let $E$ be a generic complete intersection of two surfaces of degrees $a+1, a+2$, containing
$D$. The image of the surface of degree $a+1$ is a minimal generator in $I_Y$. Let $C$ be the residual
curve to
$D$ in
$E$. By Lemma~\ref{smooth}, we can assume that both surfaces are smooth and integral. Moreover, the linear
system of curves that we obtain fixing one of the surfaces and letting the other one vary is
basepoint-free. Then
$C$ is smooth and integral by Bertini's Theorem; $C$ is non-aCM since it's CI-linked to
$D$ non-aCM. 

Applying Proposition 5.2.10 in \cite{M} to the general sections $X,Y$ of
$C,D$, we have that the minimal free resolution of $X$ is
$$0\lra R(-2a-1)\oplus R(-a-2)\lra R(-a-1)^3\lra I_X\lra 0,$$
so $X$ has degree matrix $M$.
\end{proof}

We are now ready to prove the main result of this section.
We are going to construct an integral, smooth, non-aCM curve $C\subset\PP^3$ whose general plane
section has degree matrix $M$, for any degree matrix $M$ of size at least $3\times 4$ that
has at least one entry smaller than or equal to $2$.

\begin{thm}\label{main-smooth}
Let $M=(a_{i,j})$ be an integral degree matrix, of size 
$t\times (t+1)$ such that $a_{t,1}\leq 2$, $t\geq 3$. Then
there exists an integral, smooth, non-aCM curve 
$C\subset\PP^3$ whose general plane section $X\subset\PP^2$ has degree 
matrix $M$.
\end{thm}

\begin{proof}
For some $(k,l)$, $a_{k,l}\leq 2$ and $a_{k,l+1}>0, a_{k-1,l}>0$. Fix one of such pairs $(k,l)$, and 
assume that $1\leq l\leq k-2$ and $3\leq k\leq t$. Notice that we can find such a pair $(k,l)$, since 
$M$ has positive subdiagonal by assumption.

Let $N$ be the transpose about the anti-diagonal of the first $t-1$ columns of $M$,
$$N=\left( \begin{array}{cccc} 
a_{t,t-1} & \cdots & \cdots & a_{1,t-1} \\
\vdots & & & \vdots \\
a_{t,1} & \cdots & \cdots & a_{1,1} 
\end{array} \right).$$
$N$ is a degree matrix, since $a_{2,1},\ldots,a_{t,t-1}>0$ by assumption. Moreover $N$ has one entry smaller
than $3$ since $a_{t,1}\leq 2$. 
Let $D$ be the curve constructed as in Theorem~\ref{mmain} and in particular as seen in Remark~\ref{nice},
starting from the submatrix 
$$L=\left(\begin{array}{ccc} a_{k,l+1} & a_{k-1,l+1} & a_{k-2,l+1} \\
a_{k,l} & a_{k-1,l} & a_{k-2,l}\end{array}\right).$$
Here $l$ and $k$ are the pair of integers chosen in the beginning. $L$ is a submatrix of $N$, since 
$l+1\leq t-1$ and $k-2\geq 1$.

$D$ is a non-degenerate, reduced curve with two connected components, and
it has singularities only in the intersections of its irreducible components, as we saw in
Remark~\ref{singlocusC}. It is non arithmetically Cohen-Macaulay, and its general plane section has 
degree matrix $N$.

In case $N$ is a $2\times 3$ matrix whose entries are all equal to $1$, we can still let $D$ be the
generic union of a line and a smooth plane conic. The general plane section of $D$ consists of 
three non-collinear points, hence it has degree matrix $N$. In this case, $D$ is non-degenerate, smooth, 
disconnected and non arithmetically Cohen-Macaulay. Its saturated ideal is generated in degree $2$.

The highest degree of a minimal generator of the ideal of $D$ is 
$a_{t-1,t-1}+\ldots +a_{1,1}+a_{t,k}-a_{k,k}+1$, as we showed in Remark~\ref{nice}. 
Since $a_{t,k}\leq a_{k,k}$ and $1\leq a_{t,t}$, then 
$a_{1,1}+\ldots +a_{t,t}\geq a_{t-1,t-1}+\ldots +a_{1,1}+a_{t,k}-a_{k,k}+1$. 
From Remark~\ref{smoothnice}, there exists a smooth surface $U$ of degree $a_{1,1}+\ldots+a_{t,t}$ 
containing $D$. Let $T$ be a generic surface of degree $a_{1,1}+\ldots a_{t-1,t-1}+a_{t,t+1}$. 
Abusing notation,  we refer to both the surface and its equation by $U$, or $T$. Then $I_E=(U,T)$ is the 
saturated ideal of a complete intersection $E$, containing $D$. Let $C$ be the residual curve to $D$ in $E$. 
By Bertini's Theorem, $C$ is smooth and connected. In fact, it is the general element
of the linear system of curves cut out on the smooth surface $U$, outside of $D$, by
surfaces of degree $a_{1,1}+\ldots a_{t-1,t-1}+a_{t,t+1}$. The linear system is basepoint-free, since 
$$a_{1,1}+\ldots +a_{t-1,t-1}+a_{t,t+1}\geq a_{1,1}+\ldots +a_{t-1,t-1}+a_{t,k}-a_{k,k}+1$$ 
that is bigger than or equal to the highest degree of a minimal generator of $I_D$.
The following Claim concludes the proof.

{\bf Claim:} $M$ is the degree matrix of the general plane section of $C$.

Let $X\subset\PP^2$ be the general plane section of $C$. By construction, $X$ is CI-linked to the general 
plane section $Y$ of $D$ via a $CI(a_{1,1}+ \ldots +a_{t,t},a_{1,1}+\ldots a_{t-1,t-1}+a_{t,t+1})$.
The minimal free resolution of $I_Y$ is
$$0 \ra \bigoplus_{i=1}^{t-1} R(-\sum_{j=1}^{t-1} a_{t-j,t-j} -a_{t,i}) \ra
\bigoplus_{i=0}^{t-1} R(-\sum_{j=1}^i a_{t+1-j,t-j}-\sum_{j=i+1}^{t-1} a_{t-j,t-j}) 
\ra I_Y \ra 0.$$ By Proposition 5.2.10 in \cite{M}, the minimal free resolution of $I_X$ is of the form
$$0 \lra \bigoplus_{i=1}^{t-1} R(-\sum_{j=1}^{t} a_{j,j} -a_{t,i}) \lra
\bigoplus_{i=0}^{t} R(-\sum_{j=1}^i a_{j,j}-\sum_{j=i+1}^{t} a_{j,j+1}) 
\lra I_X \lra 0.$$
This proves that the degree matrix of $X$ is $M$: no cancellation can occur in the free resolution of $X$. 
In fact, no cancellation occurs between the shifts corresponding to the submatrix $N$. The entries in 
the last two columns of $M$ are positive, since $a_{t,t}>0$, therefore no cancellation can occur there 
either.
\end{proof}

The $h$-vectors of zero-dimensional schemes of $\PP^2$ that occur as the general plane section of some
integral, smooth, curve $C\subset\PP^3$ have been characterized in \cite{GP},
\cite{Sa}, \cite{MR}, and \cite{GM1}. They are the ones of {\em decreasing type}, i.e. the $h$-vectors
$h(z)=1+h_1 z+\ldots+h_s z^s$, $h_s\neq 0$, for which
$h_i>h_{i+1}$ implies $h_{i+1}>h_{i+2}$, for $i\leq s-2$.
The results we mentioned, together with Corollary~\ref{h-vect}, Theorem~\ref{smooth2x3} and
Theorem~\ref{main-smooth}, imply the following result. 

\begin{cor}\label{hvects}
Let $h(z)=1+h_1 z+\ldots+h_s z^s$, $h_s\neq 0$, be the $h$-vector of some zero-dimensional scheme $X\subset\PP^2$.
$h(z)$ is the $h$-vector of the general plane section of some integral, smooth, non-aCM curve
$C\subset\PP^3$ if and only if it is of decreasing type and it is different from the $h$-vector of a
$CI(a,b)$, $a\neq 2$, $b\geq a$, from the $h$-vectors of Corollary~\ref{h-vect}, and from $1+2z$.
\end{cor}

\begin{proof}
If $h(z)$ is the $h$-vector of the general plane section of some integral, smooth, non-aCM curve
$C\subset\PP^3$, then it is of decreasing type, as shown in \cite{H}. Moreover, it is different from the
$h$-vector of a $CI(a,b)$, $a\neq 2$, $b\geq a$ and from the $h$-vectors of Corollary~\ref{h-vect}. In
fact, a zero-dimensional scheme that has the $h$-vector of a $CI(a,b)$ is a $CI(a,b)$, and if $a\neq 2$, $2\neq b\geq
a$ then $C$ is aCM by Theorem~\ref{strano}. If the general plane section of an integral $C$ is a
$CI(1,2)$, then $C$ is aCM. If the general plane section of $C$ has one of the $h$-vectors of
Corollary~\ref{h-vect}, then $C$ has to be arithmetically Cohen-Macaulay, by Corollary~\ref{h-vect}.

Conversely, let $h(z)$ be an $h$-vector of decreasing type, different from the $h$-vector of a
$CI(a,b)$, $a\neq 2$, $b\geq a$ and from the $h$-vectors of Corollary~\ref{h-vect}.
To any $h$-vector $h(z)$, we can uniquely associate a degree matrix $M$ with no entries equal to $0$, such
that if $X\subset\PP^2$ is a zero-dimensional scheme with degree matrix $M$, then the $h$-vector of $X$ is $h(z)$. 
Under our assumptions, $M$ can be any one of the following: \begin{itemize}
\item $M=(2,a)$ for some $a\geq 2$,
\item $M$ is a matrix of size $2\times 3$, with positive entries (since $M$ is the degree matrix of
points in Uniform Position), that do not satisfy the hypothesis of either Proposition~\ref{aCM2x3}, or
Proposition~\ref{aCM-II}, and not all of its entries are equal to $1$ (since $h(z)\neq 1+2z$),
\item $M$ is integral and has size $t\times (t+1)$, for some $t\geq 3$.
\end{itemize}
See the definition of integral degree matrix before Theorem~\ref{htv}.

If $M=(2,a)$ for some $a\geq 2$, let $C$ be a generic, smooth, rational curve on a smooth quartic
surface, as in Example~\ref{rat} or Remark~\ref{rat2}. The general plane section of $C$ has degree matrix
$M$, hence $h$-vector $h(z)$.
If $M$ is a degree matrix of size $2\times 3$ with positive entries, such that
$a_{2,1}\leq 2$, then by Theorem~\ref{smooth2x3} there exists a smooth, integral, non-aCM curve $C$
whose general plane section has degree matrix $M$, hence $h$-vector $h(z)$.
If $M$ has size bigger than or equal to $3\times 4$, and $a_{t,1}\leq 2$, then by
Theorem~\ref{main-smooth} there exists a smooth, integral, non-aCM curve $C$ whose general plane section
has degree matrix $M$, hence
$h$-vector $h(z)$.
If $M=(a_{i,j})_{i=1,\ldots,t;\; j=1,\ldots, t+1}$ has $a_{t,1}\geq 3$, let
$N=(b_{i,j})_{i=1,\ldots,t+1;\; j=1,\ldots,t+2}$ be the degree matrix with entries $b_{i,j}=a_{i,j-1}$ for
$i=1,\ldots,t$, $j=2,\ldots,t+2$, $b_{t+1,1}=0$, $b_{t+1,2}=2$. $N$ is determined by these entries, under
the assumption that it's homogeneous. $b_{i,j}>0$ for $(i,j)\neq (t+1,1)$, so $N$ is an integral degree
matrix. Moreover, the $h$-vector of a zero-dimensional scheme that has degree matrix $N$ is again $h(z)$. Then, by 
Theorem~\ref{main-smooth}, there exists a non-aCM, reduced, connected curve $C\subset\PP^3$ whose 
general plane section $X\subset\PP^2$ has degree matrix $M$, hence $h$-vector $h(z)$.  
\end{proof}

\section{Arithmetically Buchsbaum curves}

In this section we work over a field $k$ of arbitrary characteristic. We will examine the case of
arithmetically Buchsbaum curves. In particular, we will address the question of which graded Betti
numbers can correspond to points that are the general hyperplane section of an arithmetically Buchsbaum
curve. We will give an explicit characterization of the degree matrices that correspond to such points,
in the case of space curves and points in the plane. In the case of points in $\PP^n$, we will
find a necessary condition on the lifting matrix. Moreover, we will prove some bounds on the dimension of
the deficiency module
$\M_C$ of a Buchsbaum curve
$C$ and on the initial and final degree of $\M_C$, in terms of the entries of the lifting matrix of the
general plane section $X$ of $C$. We will then prove that the bounds are sharp.

\begin{defn}
Let $C\subset \PP^{n+1}$ be a curve. $C$ is \emph{arithmetically Buchsbaum}, or briefly {\em Buchsbaum},
if its deficiency module $\M_C$ is annihilated by the irrelevant maximal ideal $\MM=(x_0,\ldots,x_{n+1})$
of $S$, i.e. if its coordinate ring is Buchsbaum. 
\end{defn}

For an introduction to Buchsbaum curves and their properties, or Buchsbaum rings, see Chapter 3 of
\cite{M}, or \cite{SV}. For results about arithmetically Buchsbaum curves and their general hyperplane
section, especially in the case of space curves, see \cite{GM} and \cite{GM2}.

We begin with some observations about the deficiency module of a Buchsbaum curve.
For the whole section, $C$ will denote an arithmetically Buchsbaum curve in $\PP^{n+1}$ and $\M_C$ its
deficiency module. $X\subset\PP^n$ will be a general hyperplane section of $C$, by a  hyperplane with
equation $L=0$.

\begin{prop}\label{Buch}
Let $C\subset \PP^{n+1}$ be a Buchsbaum curve and $X\subset \PP^n$ its hyperplane section, by a general
hyperplane with equation $L=0$. Then $$\M_C=I_X/(I_C+(L))(1)\subseteq soc\; H^1_*(\I_X)$$ where $soc\;
H^1_*(\I_X)$ denotes the socle of $H^1_*(\I_X)$.
\end{prop}

\begin{proof}
Look at the short exact sequence of ideal sheaves $$0\lra \I_C(-1)\lra \I_C \lra \I_X\lra 0.$$
Taking global sections, we get the standard long exact sequence of cohomology modules
$$0\lra I_C(-1)\stackrel{\cdot L}{\lra} I_C \lra I_X\lra \M_C(-1) \stackrel{0}{\lra} \M_C \lra
H^1_*(\I_X) \lra \cdots.$$ The first map is multiplication by $L$. The map $\M_C(-1)\lra \M_C$ is
again multiplication by $L$, hence the zero map, since $C$ is Buchsbaum.

From the long exact sequence above, we can conclude that:
\begin{itemize} 
\item $\M_C(-1)=ker(\M_C(-1)\stackrel{0}{\lra} \M_C)=coker(I_C \lra I_X)=I_X/(I_C+(L))$
\item $\M_C = soc\; \M_C \subseteq soc\; H^1_*(\I_X)$.
\end{itemize}
Putting these facts together gives the thesis.
\end{proof}

\begin{cor}
With the notation of Proposition~\ref{Buch}, let 
$$0\lra\FF_n=\bigoplus_{i=1}^{t}R(-m_i)\lra\FF_{n-1}\lra\cdots\lra\FF_2\lra\FF_1\lra I_X\lra 0$$
be the minimal free resolution of $I_X$. Then
$$dim_k\; \M_C\leq t.$$
\end{cor}

Consider a zero-dimensional scheme $X\subset\PP^n$ that is a general hyperplane section of an 
arithmetically Buchsbaum, non-aCM curve $C\subset\PP^{n+1}$.
We now prove a necessary condition on the entries of the lifting matrix of $X$.

\begin{prop}\label{BuchP^n}
Let $X\subset\PP^n$ be a general hyperplane section of an arithmetically Buchsbaum, non-aCM curve
$C\subset\PP^{n+1}$. Let $M=(a_{ij})_{i=1,\ldots,t;\; j=1,\ldots,r}$ be the lifting matrix of $X$.
Then $a_{ij}=n$, for some $i,j$.
\end{prop}

\begin{proof}
By Proposition~\ref{Buch}, the deficiency module $\M_C$ of $C$ is $$\M_C=I_X/(I_C+(L)/(L))(1)\subseteq
soc\; H^1_*(\I_X)$$ where $L$ is a general linear form and $soc\; H^1_*(\I_X)$ denotes the socle of the
module
$H^1_*(\I_X)$. Since $C$ is non-aCM, $\M_C\neq 0$. $$\M_C=I_X/(I_C+(L)/(L))(1),$$ therefore 
$\alpha(\M_C)=d_j-1$ for some
$j=1,\ldots,r$. Moreover,
$$M_C\subseteq soc\; H^1_*(\I_X)=\bigoplus_{i=1}^t k(-m_i+n+1),$$
so $\alpha(\M_C)=m_i-n-1$ for some $i=1,\ldots,t$. Then $d_j=m_i-n$ for some $i,j$.
\end{proof}

We quote a result of A.V. Geramita and J. Migliore that gives a bound on the degrees of a minimal
generating system for $C\subset\PP^3$, in terms of the degrees of the minimal generators of the saturated
ideal of the general plane section $X$. We will use this result in the proof of the next theorem.

\begin{prop}(Corollary 2.5, \cite{GM2})\label{gmm}
Let $C\subset\PP^3$ be an arithmetically Buchsbaum curve, and let $X\subset\PP^2$ be its general plane
section. If $I_X$ is generated in degree less than or equal to $d$, then $I_C$ is generated in degree
less than or equal to $d+1$.
\end{prop}

We can now give a characterization of the matrices with integer
entries that occur as degree matrix of the general plane section of an arithmetically Buchsbaum, non-aCM 
curve $C\subset\PP^3$. This is a refinement of Proposition~\ref{BuchP^n}, since if $X\subseteq\PP^2$ 
then its lifting matrix and degree matrix coincide (see
Definition~\ref{liftingmatrix} and the following observations). 

\begin{thm}\label{mainBuch} 
Let $M=(a_{i,j})_{i=1,\ldots,t;\; j=1,\ldots t+1}$ be a degree matrix. Then $M$ is the degree matrix of 
the general plane section of an arithmetically Buchsbaum, non arithmetically Cohen-Macaulay curve 
$C\subset\PP^3$ if and only if $a_{i,j}=2$, for some $i,j$. 

For any such $M$, $C$ can be chosen 
such that if the ideal $I_C$ is minimally generated in degree less than or equal to $d$, then $C$ lies on 
a smooth surface of degree $d$.
\end{thm}
\begin{proof}
Assume that $M=(a_{i,j})$ is the degree matrix of some zero-dimensional scheme $X\subset\PP^2$ that is the general 
plane section of an arithmetically Buchsbaum curve $C\subset\PP^3$. Proposition~\ref{BuchP^n} proves that
$a_{i,j}=2$ for some $i,j$.

Conversely, we are going to show that the condition $a_{i,j}=2$ for some $i,j$ is sufficient in order for
$M=(a_{i,j})$ to occur as the degree matrix of the general plane section of some arithmetically Buchsbaum,
non arithmetically Cohen-Macaulay
curve $C\subset\PP^3$. We proceed by induction on the size $t$ of $M$. For each $M$, we are going to
construct a curve in the linkage class of two skew lines, i.e. a curve whose deficiency module is
one-dimensional as a $k$-vector space. 

If $t=1$, then either $M=(1,2)$ or $M=(2,a)$ for $a\geq 2$.
If $M=(1,2)$, let $C$ be two skew lines: its general plane section consists of two distinct points,
hence a $CI(1,2)$ as desired. $C$ lies on a smooth quadric surface. Since $S$ is smooth and its ideal is
generated in degree $2$, by Lemma~\ref{smooth} it lies on a smooth surface of degree $d$ for any $d\geq 2$.
If $M=(2,a)$, let $D$ consist of two skew lines, $D\subset CI(2,a+1)$. We can let the surface of degree $2$ 
be smooth, and the surface of degree $a+1$ generic. Let $C$ be the
residual curve to $D$ in the link. By Bertini's Theorem, $C$ is smooth and connected. 
Moreover, the general plane section $X$ of $C$ is linked to a $CI(1,2)$
via a $CI(2,a+1)$. Using Proposition 5.2.10 in \cite{M}, the minimal free resolution of $X$ is 
$$0\lra R(-a-2)\lra R(-a-1)\oplus R(-2)\lra I_X\lra 0$$ so $X$ is a $CI(2,a)$.
Notice that, in this case, $\M_C=k(1-a)$ and the module lies in the highest degree possible for a fixed $a$.
The ideal $I_C$ is generated in degree less than or equal to $a+1$, so $C$ lies on a smooth surface of 
degree $d$ for any $d\geq a+1$ by Lemma~\ref{smooth}.

Let $M=(a_{i,j})_{i=1,\ldots,t;\; j=1,\ldots, t+1}$ and assume that $a_{i,j}=2$ for some $2\leq j\leq t$.
Let $$N=\left( \begin{array}{cccc} 
a_{t,t} & \cdots & \cdots & a_{1,t} \\
\vdots & & & \vdots \\
a_{t,2} & \cdots & \cdots & a_{1,2} 
\end{array} \right).$$
$N$ is the transpose about the anti-diagonal of the submatrix obtained by deleting the
first and last columns of $M$. Notice that $N$ is a degree matrix.
By the induction hypothesis, there is an arithmetically Buchsbaum curve $D$ in the linkage class of
two skew lines, whose general plane section $Y$ has degree matrix $N$. The saturated ideal
$I_D$ of $D$ is generated in degree less than or equal to $a_{1,2}+\ldots+a_{t-1,t}+1$, by
Proposition~\ref{gmm}.\newline
$\alpha(I_D)\leq a_{2,2}+\ldots+a_{t,t}+1$. So we can find a complete intersection of forms of degrees
$a_{1,1}+\ldots +a_{t,t},a_{1,t+1}+a_{2,2}+\ldots+a_{t,t}$ containing $D$.
Both the surfaces that cut out the complete intersection can be chosen in such a way that their
images in $I_Y$ are not minimal generators. Let $C$ be the residual
of $D$ in the $CI(a_{1,1}+\ldots +a_{t,t},a_{1,t+1}+a_{2,2}+\ldots+a_{t,t})$. 
By the Hartshorne-Schenzel Theorem, $C$ is in the linkage class of two skew lines, as is $D$. Since $Y$
has degree matrix $N$, using Proposition 5.2.10 in \cite{M}, we see that $X$ has degree matrix $M$.
The surface of degree $a_{1,t+1}+a_{2,2}+\ldots+a_{t,t}$ can be taken to be smooth, by the induction 
hypothesis applied to $D$. The ideal of $C$ is generated in degree less than or equal to 
$a_{1,t+1}+a_{2,2}+\ldots+a_{t,t}+1$, by Proposition~\ref{gmm}. Let $d\geq a_{1,t+1}+a_{2,2}+
\ldots+a_{t,t}+1$, and consider the linear system $\Delta_d$ of surfaces of degree $d$ containing $C$.
We want to show that the general element is smooth. By Bertini's Theorem, it is smooth outside of $C$.
Consider now a point $P\in C$. By Corollary 2.10 in \cite{GV}, it is enough to exhibit two elements of 
$\Delta_d$ meeting transversally at $P$. If $C$ is smooth at $P$, we have two minimal generators 
of $I_C$, call them $F$ and $G$, meeting transversally at $P$. The degree of each of them is at most $d$. 
Add generic planes as needed, to obtain surfaces of degree $d$ that meet transversally at $P$.
Finally, we need to check that the singular points of $C$ are not fixed singular
points for $\Delta_d$. So it is enough to find a surface for each of those points that contains $C$ and is
non-singular at $P$. This follows from the fact that we have a smooth surface containing $C$ of degree 
$a_{1,t+1}+a_{2,2}+\ldots+a_{t,t}<d$. Add generic planes as needed to get a surface that is non-singular at 
$P$ and contains $C$.

Consider now the case $a_{i,1}=2$ for some $i\neq 1$.
Let $$N=\left( \begin{array}{cccc} 
a_{t,t} & \cdots & \cdots & a_{1,t} \\
\vdots & & & \vdots \\
a_{t,3} & \cdots & \cdots & a_{1,3} \\
a_{t,1} & \cdots & \cdots & a_{1,1} 
\end{array} \right),$$ 
$N$ is the transpose about the anti-diagonal of the matrix obtained deleting the second and last column
of $M$. Notice that $N$ is a degree matrix, since $a_{2,1}\geq a_{i,1}>0$. 
By the induction hypothesis, there is an arithmetically Buchsbaum curve $D$ in the
linkage class of two skew lines, whose general plane section $Y$ has degree matrix $N$. The saturated
ideal $I_D$ of $D$ is generated in degree less then or equal to $a_{1,1}+a_{2,3}+\ldots+a_{t-1,t}+1$, by
Proposition~\ref{gmm}. $\alpha(I_D)\leq a_{t,t}+\ldots+a_{3,3}+a_{2,1}+1$, so we can find
a complete intersection of forms of degrees
$a_{t,t}+\ldots +a_{3,3}+a_{2,1}+a_{1,2}, a_{t,t}+\ldots +a_{3,3}+a_{2,1}+a_{1,t+1}$ containing $D$. 
Both the surfaces that cut out the complete intersection can be chosen in such a way that their
images in $I_Y$ are not minimal generators. Let $C$ be the residual of $D$ in the
$CI(a_{t,t}+\ldots +a_{3,3}+a_{2,1}+a_{1,2}, a_{t,t}+\ldots +a_{3,3}+a_{2,1}+a_{1,t+1})$.
$C$ is in the linkage class of two skew lines, by the Hartshorne-Schenzel Theorem.
Since $Y$ has degree
matrix $N$, using Proposition 5.2.10 in \cite{M}, we see that $X$ has degree matrix $M$.
The surface of degree $a_{t,t}+\ldots +a_{3,3}+a_{2,1}+a_{1,t+1}$ can be taken smooth, by induction 
hypothesis applied to $D$. The ideal of $C$ is generated in degree less than or equal to 
$a_{1,t+1}+a_{2,2}+\ldots+a_{t,t}+1$, by Proposition~\ref{gmm}. Let 
$d\geq a_{1,t+1}+a_{2,2}+\ldots+a_{t,t}+1$, and consider the linear system $\Delta_d$ of surfaces of 
degree $d$ containing $C$.
We want to show that the general element is smooth. By Bertini's Theorem, it is smooth outside of $C$.
Consider now a point $P\in C$. By Corollary 2.10 in \cite{GV}, it is enough to exhibit two elements of 
$\Delta_d$ meeting transversally at $P$. If $C$ is smooth at $P$, we have two minimal generators 
of $I_C$, call them $F$ and $G$, meeting transversally at $P$. The degree of each of them is at most $d$. 
Add generic planes as needed, to obtain surfaces of degree $d$ that meet transversally at $P$.
Finally, we need to check that the singular points of $C$ are not fixed singular
points for $\Delta_d$. So it is enough to find a surface for each of those points that contains $C$ and is
non-singular at $P$. This follows from the fact that we have a smooth surface containing $C$ of degree 
$a_{t,t}+\ldots +a_{3,3}+a_{2,1}+a_{1,t+1}\leq d$. 
Add generic planes as needed to get a surface that is non-singular at $P$ and contains $C$.

Assume now that $a_{1,1}=2$, i.e. $i=j=1$. Let 
$$N=\left( \begin{array}{cccc} 
a_{1,1} & \cdots & \cdots & a_{1,t} \\
\vdots & & & \vdots \\
a_{t-1,1} & \cdots & \cdots & a_{t-1,t}
\end{array} \right),$$ 
be the submatrix of $M$, consisting of the first $t-1$ rows and first $t$ columns. 
By the induction hypothesis, there is an arithmetically Buchsbaum curve $D$ in the
linkage class of two skew lines, whose general plane section $Y$ has degree matrix $N$. The saturated
ideal $I_D$ of $D$ is generated in degree less than or equal to $a_{1,2}+\ldots+a_{t-1,t}+1$, by
Proposition~\ref{gmm}.
$\alpha(I_D)\leq a_{1,1}+\ldots+a_{t-1,t-1}+1$, so we can find a surface $S$ of degree
$s=a_{1,1}+\ldots+a_{t,t}$, containing $D$. The surface can be chosen such that its image in $I_Y$ is not 
a minimal generator.
Perform a basic double link of degrees $s,a_{t,t+1}$.
Let $C$ be the curve obtained in the $BDL(a_{1,1}+\ldots+a_{t,t},a_{t,t+1})$. Let the surface $F$ of 
degree $a_{t,t+1}$ be generic.
$C$ is in the linkage class of two skew lines, as $D$ is. Since $Y$ has degree matrix $N$, using
Proposition 5.4.5 in \cite{M}, we see that $X$ has degree matrix $M$. No cancellation can occur, since
the image of $S$ in $I_Y$ is not a minimal generator, and by genericity of $F$.
The ideal of $C$ is generated in degree less than or equal to 
$a_{1,t+1}+a_{2,2}+\ldots+a_{t,t}+1$, by Proposition~\ref{gmm}. Let 
$d\geq a_{1,t+1}+a_{2,2}+\ldots+a_{t,t}+1$, and consider the linear system $\Delta_d$ of surfaces of 
degree $d$ containing $C$.
We want to show that the general element is smooth. By Bertini's Theorem, it is smooth outside of $C$.
Consider now a point $P\in C$. By Corollary 2.10 in \cite{GV}, it is enough to exhibit two elements of 
$\Delta_d$ meeting transversally at $P$. If $C$ is smooth at $P$, we have two minimal generators 
of $I_C$, call them $F$ and $G$, meeting transversally at $P$. The degree of each of them is at most $d$. 
Add generic planes as needed, to obtain surfaces of degree $d$ that meet transversally at $P$.
Finally, we need to check that the singular points of $C$ are not fixed singular
points for $\Delta_d$. So it is enough to find a surface for each of those points that contains $C$ and is
non-singular at $P$. By the induction hypothesis, we can find a smooth surface $T$ of degree 
$a_{1,1}+a_{2,3}+\ldots+a_{t-1,t}+1$ containing $D$. By genericity, we can assume that the surface $F$
used in the construction of $C$ is smooth. $T\cup F$ is a surface of degree 
$a_{1,1}+a_{2,3}+\ldots+a_{t,t+1}+1=a_{1,t+1}+a_{2,1}+a_{3,3}+\ldots+a_{t,t}<d.$
Add generic planes as needed to get a surface that is non-singular at each point of $C$, except for the 
points of intersection of $D$ and $S\cap F$. The surfaces $S$ and $T\cup F$ meet transversally, so those 
cannot be fixed singular points of $\Delta_d$ either.

Finally, let $j=t+1$, i.e. $a_{i,t+1}=2$ for some $i$. Let
$$N=\left( \begin{array}{cccc} 
a_{t,t+1} & \cdots & \cdots & a_{1,t+1} \\
a_{t,t-1} & \cdots & \cdots & a_{1,t-1} \\
\vdots & & & \vdots \\
a_{t,2} & \cdots & \cdots & a_{1,2}
\end{array} \right),$$
$N$ is the transpose about the anti-diagonal of the matrix obtained deleting the first and $t$-th columns 
of $M$. By the induction hypothesis, there is an arithmetically Buchsbaum curve $D$ in the
linkage class of two skew lines, whose general plane section $Y$ has degree matrix $N$. The
saturated ideal $I_D$ of $D$ is generated in degree less than or equal to
$a_{1,2}+\ldots+a_{t-2,t-1}+a_{t-1,t+1}+1$, by Proposition~\ref{gmm}. Moreover,\newline 
$\alpha(I_D)\leq a_{t,t+1}+a_{t-1,t-1}+\ldots+a_{2,2}+1$, so we can find
a complete intersection of forms of degrees
$a_{t,t+1}+a_{t-1,t-1}+\ldots +a_{1,1}, a_{t,t+1}+a_{t-1,t-1}+\ldots +a_{2,2}+a_{1,t}$ containing $D$.
Both the surfaces that cut out the complete intersection can be chosen in such a way that their
images in $I_Y$ are not minimal generators. Let $C$ be the residual curve to $D$ in the 
$CI(a_{t,t+1}+a_{t-1,t-1}+\ldots +a_{1,1}, a_{t,t+1}+a_{t-1,t-1}+\ldots +a_{2,2}+a_{1,t})$.
By the Hartshorne-Schenzel Theorem, $C$ is in the linkage class of two skew lines, as is $D$.
Since $Y$ has degree matrix $N$, using Proposition 5.2.10 in \cite{M}, we see that $X$ has degree 
matrix $M$.
The surface of degree $a_{t,t+1}+a_{t-1,t-1}+\ldots +a_{2,2}+a_{1,t}$ can be taken to be smooth, by the
induction hypothesis applied to $D$. The ideal of $C$ is generated in degree less than or equal to 
$a_{1,t+1}+a_{2,2}+\ldots+a_{t,t}+1$, by Proposition~\ref{gmm}. Let 
$d\geq a_{1,t+1}+a_{2,2}+\ldots+a_{t,t}+1$, and consider the linear system $\Delta_d$ of surfaces of 
degree $d$ containing $C$.
We want to show that the general element is smooth. By Bertini's Theorem, it is smooth outside of $C$.
Consider now a point $P\in C$. By Corollary 2.10 in \cite{GV}, it is enough to exhibit two elements of 
$\Delta_d$ meeting transversally at $P$. If $C$ is smooth at $P$, we have two minimal generators 
of $I_C$, call them $F$ and $G$, meeting transversally at $P$. The degree of each of them is at most $d$. 
Add generic planes as needed, to obtain surfaces of degree $d$ that meet transversally at $P$.
Finally, we need to check that the singular points of $C$ are not fixed singular
points for $\Delta_d$. So it is enough to find a surface for each of those points that contains $C$ and is
non-singular at $P$. This follows from the fact that we have a smooth surface containing $C$ of degree 
$a_{t,t+1}+a_{t-1,t-1}+\ldots +a_{2,2}+a_{1,t}<d$. 
Add generic planes as needed to get a surface that is non-singular at $P$ and contains $C$.
\end{proof}

Let us observe a few consequences of the theorem we just proved.

\begin{rem}
The proof of Theorem~\ref{mainBuch} shows that the following facts about a degree matrix
$M=(a_{i,j})_{i=1,\ldots,t;\; j=1,\ldots t+1}$ are equivalent:
\begin{itemize}
\item $a_{i,j}=2$ for some $i,j$;
\item there exists a zero-dimensional scheme $X\subset \PP^2$ and a Buchsbaum, non-aCM curve $C\subset\PP^3$ such that
$X$ is the general plane section of $C$ and $M$ is the degree matrix of $X$;
\item there exists a zero-dimensional scheme $X\subset \PP^2$ and a Buchsbaum curve $C\subset\PP^3$ in the linkage
class of two skew lines such that $X$ is the general plane section of $C$ and $M$ is the degree matrix of
$X$.
\end{itemize}
\end{rem}

\begin{rem}
Introducing a minor modification in the proof, we can show that we can always construct a curve $C$ whose
deficiency module is $\M_C=k(-d_m+1)$, where $m=\mmin\{j\; | \;
a_{i,j}=2,\; \mbox{for some $i$}\}$. Notice that this is the highest possible degree in which the deficiency 
module can appear, for a given degree matrix (see also Proposition~\ref{BuchBounds}).
\end{rem}

\begin{rem}
Theorem~\ref{mainBuch} also proves that $d=\bin{n}{2}$ generic points in $\PP^2$ cannot be the general
plane section of an arithmetically Buchsbaum curve for any $n$, unless the curve is arithmetically
Cohen-Macaulay. This had been observed already by A.V. Geramita and J. Migliore in \cite{GM}, Proposition
4.9.
\end{rem}

Our result extends a result by A.V. Geramita and J. Migliore for arithmetically Buchsbaum curves in
$\PP^3$. In \cite{GM}, they prove the following.

\begin{prop}(\cite{GM}, Proposition 4.7)
Let $C\subset\PP^3$ be an arithmetically Buchsbaum curve lying on no quadric surface. Let $C\cap L$ be a
general plane section. Assume that $\alpha(I_C)=\alpha(I_{C\cap L})$ and that $C\cap L$ is a complete
intersection. Then $C$ is a complete intersection.
\end{prop}

We are now going to consider the case of integral, arithmetically Buchsbaum curves in $\PP^3$. We want to
investigate which degree matrices can occur for a general
plane section of an integral, arithmetically Buchsbaum curve. 

\begin{notat}\label{notBuch}
Let $C\subset\PP^3$ be an integral, Buchsbaum curve, and let $X\subset\PP^2$ be its general plane section.
Let $M=(a_{i,j})_{i=1,\ldots,t;\; j=1,\ldots,t+1}$ be the degree matrix of $X$. $M$ is then an integral
matrix, i.e. $a_{i+1,i}>0$ for all $i$ (see the introduction of Section 3). By Theorem~\ref{mainBuch},
$a_{i,j}=2$ for some $i,j$.
\end{notat}

\begin{rem}
In Section 3, we saw some classes of degree matrices $M$ of size $2\times 3$ such that, if the general
plane section of an integral curve $C\subset\PP^3$ has degree matrix $M$, then $C$ is forced to be aCM (see
Propositions~\ref{aCM2x3} and \ref{aCM-II}). Notice that all of those matrices have no entry equal to
$2$, so they cannot correspond to the general plane section of an arithmetically Buchsbaum curve.
\end{rem}

We can give a characterization of the matrices $M$ that occur as the degree matrix of the general plane
section of an arithmetically Buchsbaum, integral curve $C\subset\PP^3$. They have to satisfy the 
conditions of Notation~\ref{notBuch}. We are going to show that for each of these matrices the curve 
$C$ can be taken to be smooth and connected.
  
We treat separately the case $t=1$, when $X$ is a complete intersection.

\begin{rem}
Any integral curve $C\subset\PP^3$ of degree $2$ is a plane conic. So there cannot be any arithmetically
Buchsbaum curve that is non-aCM and whose general plane section is a $CI(1,2)$.
\end{rem}

\begin{prop}
Assume that $char(k)=0$.
Let $M=(a,b)$, $b\geq a>0$. $M$ is the degree matrix of the general plane section of some smooth,
integral, arithmetically Buchsbaum, non-aCM curve $C\subset\PP^3$ if and only if $a=2$.
\end{prop}

\begin{proof}
Assume that  $M$ is the degree matrix of the general plane section of some smooth,
integral, Buchsbaum, non-aCM curve $C\subset\PP^3$. We already saw that $M$ needs to contain a $2$ (see
Proposition~\ref{BuchP^n}). The Remark above shows that $a\neq 1$, so $a=2$.

Conversely, let $M=(2,b)$, $b\geq 2$. We want to construct an integral, smooth, arithmetically Buchsbaum,
non-aCM curve
$C$, whose general plane section has degree matrix $M$.
Let $D$ be two skew lines, and let $Q$ be a smooth, integral quadric surface containing $D$. Consider the
linear system of curves cut out on $Q$, outside of $D$, by surfaces of degree $b+1$ containing $D$. It is
basepoint-free, since $I_D$ is generated in degree $2<b+1$. By Bertini's Theorem (see Theorem~\ref{bert}),
the general element $C$ of the linear system is smooth and integral. $C$ is in the linkage class of two
skew lines by construction, and its general plane section has degree matrix $M$, by Proposition 5.2.10 in
\cite{M}.
\end{proof}

We are now going to characterize the integral matrices that can occur as the degree matrix of the general
plane section of an arithmetically Buchsbaum, non-aCM, integral curve of $\PP^3$. 

\begin{thm}\label{intBuch}
Assume that $char(k)=0$.
Let $M=(a_{i,j})_{i=1,\ldots,t;\; j=1,\ldots t+1}$, $t\geq 2$ be an integral degree matrix. Then $M$ is
the degree matrix of the general plane section of an
arithmetically Buchsbaum, non-aCM, integral curve $C\subset\PP^3$ if
and only if $a_{i,j}=2$, for some $i,j$.
Moreover, for such a matrix $M$ the curve $C$ can be chosen to be smooth and integral.
\end{thm}

\begin{proof}
The necessity of the hypothesis $a_{i,j}=2$ has been proven in Theorem~\ref{mainBuch}.

Let $M$ be an integral degree matrix such that $a_{i,j}=2$, for some $i,j$. We are going to construct an
integral, smooth, Buchsbaum curve $C$ in the linkage class of two skew lines, such that its general plane
section $X$ has degree matrix $M$.
We start from degree matrices of size $2\times 3$. Notice that in this case, all the entries of $M$ are
positive. Then we have the following possibilities for $M$.

\textbf{Case 1.} Let $$M=\left(\begin{array}{ccc} 2 & a & b \\
1 & a-1 & b-1 \end{array}\right)$$
and let $D$ be two skew lines. $D\subset CI(a+1,b+1)$, where the surface of degree $a+1\geq 3$ can be
taken smooth and integral. Choosing a generic surface of degree $b+1$, we have that the residual to $D$
in the complete intersection is smooth and integral by Bertini's Theorem, since the
linear system of curves cut out outside of $D$ by surfaces of degree $b+1$ containing $D$ is
basepoint-free. Let $C$ be the residual curve to $D$ in the complete intersection. 
The general plane section of $C$ has degree matrix $M$ by Proposition 5.2.10 in \cite{M}.

\textbf{Case 2.} Let $$M=\left(\begin{array}{ccc} a & b & c \\
2 & b+2-a & c+2-a \end{array}\right)$$
and let $D$ be the residual to two skew lines in a general $CI(2,a+1)$ (see the proof of
Theorem~\ref{mainBuch}). The ideal of $D$ is generated in degree less than or equal to $a+1$ and its
general plane section is a $CI(2,a)$.
$D\subset CI(b+2,c+2)$, where the surface of degree
$b+2\geq a+2$ can be taken to be smooth and integral. Choosing a generic surface of degree $c+2$, we have 
that the residual to $D$ in the complete intersection is smooth and integral by Bertini's Theorem, since 
the linear system of curves cut out outside of $D$ by surfaces of degree
$c+2$ containing $D$ is basepoint-free (because the ideal $I_D$ is generated in degree less than or equal
to $a+1$). Let $C$ be the residual curve to $D$ in the complete intersection. 
The general plane section of $C$ has degree matrix $M$ by Proposition 5.2.10 in \cite{M}.

\textbf{Case 3.} Let $$M=\left(\begin{array}{ccc} 1 & 2 & a \\
1 & 2 & a \end{array}\right)$$
and let $D$ be two skew lines.
$D$ is contained in a smooth, integral surface of degree $3$, call it $S$. Perform a basic double link on
$S$, using a general surface of degree $a$, and let $C=D\cup CI(3,a)$. The general plane
section of $C$ has degree matrix $M$ by Proposition 5.4.5 in \cite{M}.
The linear system of curves on $S$ that are linearly equivalent to $C$ is basepoint-free (in fact, the
linear system $|D|$ is itself basepoint-free, as shown in Theorem~\ref{smooth2x3}), so the general element
of
$|C|$ is smooth and integral. By Lemma~\ref{deform}, its general plane section has degree matrix $M$.

\textbf{Case 4.} Let $$M=\left(\begin{array}{ccc} 1 & 1 & 2 \\
1 & 1 & 2 \end{array}\right)$$
and let $D$ be two skew lines.
$D$ is contained in a smooth, integral surface of degree $3$, call it $S$. Perform a basic double link on
$S$, using a general plane, let $C=D\cup CI(1,3)$. The general plane
section of $C$ has degree matrix $M$ by Proposition 5.4.5 in \cite{M}.
The linear system of curves on $S$ that are linearly equivalent to $C$ is basepoint-free (in fact, the
linear system $|D|$ is itself basepoint-free, as in the proof of Theorem~\ref{smooth2x3}), so the general
element of $|C|$ is smooth and integral. By Lemma~\ref{deform}, its general plane section has degree
matrix $M$.

This concludes the proof of the case $t=2$.
Assume now that $t\geq 3$ and that $j\leq t-1$. Consider the submatrix
$$N=\left(\begin{array}{cccc} a_{t,t-1} & \cdots & \cdots & a_{1,t-1} \\
\vdots & & \vdots \\
a_{t,1} & \cdots & \cdots & a_{1,1}
\end{array}\right),$$
$N$ is the transpose about the
anti-diagonal of the first $t-1$ columns of $M$. By induction, we have an integral, smooth, Buchsbaum
curve
$D$ in the linkage class of two skew lines, whose general plane section has degree matrix $N$. The ideal
of
$D$ is generated in degree less than or equal to $a_{1,1}+\ldots+a_{t-1,t-1}+1$ (see
Proposition~\ref{gmm}). So there is a smooth surface $S$ of degree $s=a_{1,1}+\ldots+a_{t,t}$ containing
$D$, by Lemma~\ref{smooth}. Consider the linear system of curves cut out on $S$, outside of $D$, by
surfaces of degree $a_{1,1}+\ldots+a_{t-1,t-1}+a_{t,t+1}$ containing $D$. The linear system is
basepoint-free, so by Bertini's Theorem, the general element $C$ is smooth and
integral. The general plane section of $C$ has degree matrix $M$ by Proposition 5.2.10 in \cite{M}.

The cases when $j=t,t+1$ can be proved in an analogous way. If $j=t$, start from the degree
matrix
$$N=\left(\begin{array}{cccc} a_{t,t} & \cdots & \cdots & a_{1,t} \\
a_{t,t-2} & \cdots & \cdots & a_{1,t-2} \\
\vdots & & \vdots \\
a_{t,1} & \cdots & \cdots & a_{1,1}
\end{array}\right),$$
$N$ is the transpose about the
anti-diagonal of the submatrix of $M$ obtained by deleting columns $t-1$ and $t+1$. Link via a 
$CI(a_{1,1}+\ldots+a_{t,t},a_{1,1}+\ldots+a_{t-1,t-1}+a_{t,t+1})$.

If $j=t+1$, start from the degree
matrix
$$N=\left(\begin{array}{cccc} a_{t,t+1} & \cdots & \cdots & a_{1,t+1} \\
a_{t,t-2} & \cdots & \cdots & a_{1,t-2} \\
\vdots & & \vdots \\
a_{t,1} & \cdots & \cdots & a_{1,1}
\end{array}\right),$$
$N$ is the transpose about the
anti-diagonal of the submatrix of $M$ obtained by deleting columns $t-1$ and $t$. Link via a 
$CI(a_{1,1}+\ldots+a_{t-1,t-1}+a_{t,t+1},a_{1,1}+\ldots+a_{t-2,t-2}+a_{t-1,t}+a_{t,t+1})$.
\end{proof}

\begin{rem}\label{implBuch}
As in Theorem~\ref{mainBuch}, we showed that the following facts about an {\em integral} degree matrix
$M=(a_{i,j})_{i=1,\ldots,t;\; j=1,\ldots t+1}$ are equivalent:
\begin{itemize}
\item $a_{i,j}=2$ for some $i,j$;
\item there exists a zero-dimensional scheme $X\subset \PP^2$ and an integral Buchsbaum, non-aCM curve $C\subset\PP^3$
such that $X$ is the general plane section of $C$ and $M$ is the degree matrix of $X$;
\item there exist a zero-dimensional scheme $X\subset \PP^2$ and an integral, smooth Buchsbaum, non-aCM curve
$C\subset\PP^3$ such that $X$ is the general plane section of $C$ and $M$ is the degree matrix of $X$;
\item there exists a zero-dimensional scheme $X\subset \PP^2$ and an integral Buchsbaum curve $C\subset\PP^3$ in the
linkage class of two skew lines such that $X$ is the general plane section of $C$ and $M$ is the degree
matrix of $X$; 
\item there exists a zero-dimensional scheme $X\subset \PP^2$ and an integral, smooth Buchsbaum curve $C\subset\PP^3$
in the linkage class of two skew lines such that $X$ is the general plane section of $C$ and $M$ is the
degree matrix of $X$.
\end{itemize}
\end{rem}

\begin{rem}
G. Paxia and A. Ragusa proved in \cite{PR} that any integral, arithmetically Buchsbaum curve
$C\subset\PP^3$ can be deformed to a smooth integral curve. The deformation, moreover, preserves the
cohomology, hence the deficiency module, of the curve. Their proof relies heavily on papers of M.
Martin-Deschamps and D. Perrin (\cite{MDP}) and of S. Nollet (\cite{N}).

Their result is related to some of the implications of Remark~\ref{implBuch}. In fact, we show that
the existence of an integral, arithmetically Buchsbaum curve, whose general plane section has a
prescribed degree matrix is equivalent to the existence of an integral, smooth, arithmetically Buchsbaum
curve, whose general plane section has that same degree matrix. 
However, deforming an integral, arithmetically Buchsbaum curve to an integral, smooth one does not in
general preserve the degree matrix of the general plane section. In particular, the way the deformation
is done in \cite{PR} implies that if the general plane section $X$ of an integral, Buchsbaum curve $C$ has
a minimal free resolution $$0\lra\FF_2\oplus\FF\lra\FF_1\oplus\FF\lra I_X\lra 0$$ where $\FF_2$ and
$\FF_1$ are free $R$-modules without any (abstractly) isomorphic free summand, then the minimal free
resolution of the general plane section $Y$ of the smooth, integral deformation $D$ of $C$ is
$$0\lra\FF_2\lra\FF_1\lra I_Y\lra 0.$$
In particular, the result of G. Paxia and A. Ragusa does not imply the result of Theorem~\ref{intBuch}.
\end{rem}

We now turn to the study of the deficiency module of a Buchsbaum curve. The ground field $k$ can have any 
characteristic. Using Proposition~\ref{Buch},
some easy bounds for the initial and final degrees of $\M_C$ in terms of the entries of the lifting
matrix of $X$ can be derived.
From here on, we assume only that the curve $C\subset\PP^{n+1}$ is
arithmetically Buchsbaum (hence locally Cohen-Macaulay), non-aCM, equidimensional and nondegenerate. 

\begin{prop}\label{BuchBounds}
Let $C\subset\PP^{n+1},\; X\subset\PP^n$ be as above and let $M=(a_{ij})_{i=1,\ldots,t;\; j=1,\ldots r}$
be the lifting matrix of
$X$. Then $$\alpha(\M_C)\geq \mmax \{ m_t-n-1, \; \alpha(I_X)-1 \}$$ and
$$\alpha(\M_C)^+\leq m_1-n-1.$$
\end{prop}

\begin{proof}
With our notation $$soc\; H^1_*(\I_X)=\bigoplus_{i=1}^t k(-m_i+n+1).$$ So
$\alpha(\M_C)\geq m_t-n-1$ and $\alpha(\M_C)^+\leq m_1-n-1.$ Moreover,
$$\M_C=I_X/(I_C+(L))(1)$$ gives $\alpha(\M_C)\geq d_r-1=\alpha(I_X)-1.$
\end{proof}

Following the same principle, we can give a more precise estimate of what the initial degree of the
deficiency module of $C$ can be.

\begin{rem}
Since $\M_C=I_X/(I_C+(L))(1)$, then $d_m-1\leq\alpha (\M_C)\leq d_l-1$ where\newline 
$m=\mmax\{j\; | \; a_{i,j}=2,\; \mbox{for some $i$}\}$ and 
$l=\mmin\{j\; | \; a_{i,j}=2,\; \mbox{for some $i$}\}$.
\end{rem}

From Propositions~\ref{Buch} and~\ref{BuchBounds}, we can deduce a bound on the dimension of $\M_C$ in
each degree, hence a bound on the dimension of $\M_C$ as a $k$-vector space. 

\begin{prop}\label{dimMB}
Let $$J=\{\: j \; | \; d_j=m_{k(j)}-n\;\; \mbox{for some $k(j)$}\}$$
and for each $j\in J$ let $\mu(j)$ be the number of minimal generators of $I_X$ of degree $d_j$.
Then, for $i\in\ZZ$, the dimension of the $i$-th graded component of $\M_C$ is
$$dim_k(\M_C)_i=0 \;\;\; \mbox{if\: $i\neq d_j-1$ for all $j\in J$}$$
and for $j\in J$
$$dim_k(\M_C)_{d_j-1}\leq \mmin \{ dim\; soc\; H^1_*(\I_X)_{d_j-1}, \mu(j) \}.$$
Then 
$$dim_k(\M_C)\leq \sum_{j\in J} \mmin \{ dim\; soc\; H^1_*(\I_X)_{d_j-1}, \mu(j) \}.$$
\end{prop}

\begin{proof}
First we observe that the set of all degrees $i$ where we can possibly have\newline
$dim(\M_C)_i\neq 0$ is $\{d_j-1 \; | \; j\in J\}$. In fact, by Proposition~\ref{Buch} 
$$\M_C=I_X/(I_C+(L))(1)\subseteq soc\; H^1_*(\I_X).$$
In particular, $(\M_C)_i$ can be non-zero only for $i\in \{ m_1-n-1,\ldots,m_t-n-1 \}$, since those are the 
degrees in which $soc\; H^1_*(\I_X)$ is non-zero.
Clearly, each minimal generator 
of $\M_C$ is a minimal generator of $I_X/(I_C+(L))(1)$. Therefore, each minimal generator of $\M_C$ 
has degree $d_j-1$ for some $j$. Since by assumption the structure of $\M_C$ as an 
$S$-module is trivial, a minimal system of generators of $\M_C$ as an $S$-module is also a basis as 
a $k$-vector space. Then the set of all possible degrees where the deficiency module can possibly be 
non-zero is $\{d_j-1 \; | \; j\in J\}$.
Moreover, in each degree $i=d_j-1$ where $dim(\M_C)_i$ can be non-zero we have
$$dim(\M_C)_{d_j-1}\leq \mmin \{ dim\; soc\; H^1_*(\I_X)_{d_j-1}, \mu(j) \}.$$
\end{proof}

\begin{rem}
Notice that $\mu(j)$ is the number of columns that are equal to the $j$-th column. Moreover,
$dim\; soc\; H^1_*(\I_X)_{d_j-1}=dim\; soc\; H^1_*(\I_X)_{m_{k(j)}-n-1}$ is the number of rows that are 
equal to the $k(j)$-th row. Here $k(j)$ is an integer such that the entry $(j,k(j))$ of the lifting matrix 
is equal to $n$. 
\end{rem}

\begin{defn}
Let $M$ be a lifting matrix. By a {\bf block} of entries equal to $n$ we mean a group of entries
of $M$ such that:\begin{itemize}
\item $a_{i,j}=n$ for $i_1\leq i\leq i_2$ and $j_1\leq j\leq j_2$, and 
\item $a_{i,j}\neq n$ if either $i=i_1-1$ and $j_1\leq j\leq j_2$, or $i=i_2+1$ and $j_1\leq j\leq j_2$, 
or $j=j_1-1$ and $i_1\leq i\leq i_2$, or $j>j_2$ and $i_1\leq i\leq i_2$.
\end{itemize}
\end{defn}

\begin{rem}\label{blocks}
The proof of Proposition~\ref{dimMB} also shows that each block of $n$'s in the lifting matrix 
corresponds to a degree in which the deficiency module of $C$ is possibly non-zero.
\end{rem}

From our observations, we can easily derive a criterion for lifting minimal generators from the saturated 
ideal $I_X$ of a general hyperplane section $X$, to the saturated ideal $I_C$ of the curve $C$. Notice that
this sufficient condition is weaker than the sufficient condition of Lemma~\ref{lift}, for curves 
that are not necessarily Buchsbaum.

\begin{cor}\label{liftB}
Let $C$ be an arithmetically Buchsbaum, non arithmetically Cohen-Macaulay curve, let $X$ be its general 
hyperplane section, and 
let $M$ be the lifting matrix of $X$. If for some $j$ we have $a_{ij}\neq n$ for all $i$, then the 
minimal generators of degree $d_j$ of $I_X$ lift to $I_C$.
In particular, if $a_{1,j}<n$ then the minimal generators of degrees $d_1,\ldots,d_j$ of $I_X$ lift to 
$I_C$.
\end{cor}

\begin{proof}
Let $$0\lra \bigoplus_{i=1}^t R(-m_i)\lra
\ldots\lra \bigoplus_{j=1}^r R(-d_j)\lra I_X\lra 0$$
be the minimal free resolution of $I_X$. Since
$d_j=m_i-a_{ij}$, it follows that $d_j\neq m_i-n$ if and
only if $a_{ij}\neq n$. Fix a $j$ such that $a_{ij}\neq n$ for all $i$. Then $d_j\neq m_i-n$ for all $i$, so 
$(\M_C)_{d_j-1}=0$ by Proposition~\ref{dimMB}. 
Therefore all the minimal generators of degree $d_j$ of $I_X$ lift to $I_C$.
This proves the first part of the statement.

Assume now that $a_{1,j}<n$ for some $j$. Then $a_{1,l}<n$ for $l\leq j$. 
In particular, $a_{il}\neq n$ for all $i$ and for all $l\leq j$. Then the minimal generators of 
degrees $d_1,\ldots,d_j$ of $I_X$ lift to $I_C$.
\end{proof}

\begin{rem}
In the case of points in $\PP^2$,
assuming $a_{1,j}<2$ is equivalent to assuming $a_{1,j}=1$. In fact $a_{1,j}\leq 0$ implies
$a_{i,j}\leq 0$ for all $i$, and the Hilbert-Burch matrix of a scheme of codimension 2 cannot have a
column of zeroes.   
\end{rem}

\begin{rem}\label{bingen}
Corollary~\ref{liftB} clarifies how, for $X\subset\PP^2$ a generic zero-scheme of degree $d=\bin{n}{2}$
for some $n$, an arithmetically Buchsbaum curve of $\PP^3$ that has $X$ as its general plane section needs
to be arithmetically Cohen-Macaulay as well. In fact, all the entries of the degree matrix of $X$ are
equal to 1.
\end{rem}

We now look at space curves whose deficiency modules are concentrated in one degree. We see how in
this special case the bounds on the dimension of $\M_C$ of Proposition~\ref{BuchBounds} and
Corollary~\ref{dimMB} are sharp. We concentrate on minimal curves in $\PP^3$.

\begin{ex}\label{mincurves}
Let $C_n$ be a minimal curve for its Liaison class (see \cite{M} for definition and facts
about minimal curves) and let $M_{C_n}=K^n(-2n+2)$ be its deficiency module. Let $S=k[x_0,x_1,x_2,x_3]$
and
$R=k[x_0,x_1,x_2]$. We can construct such a $C_n$ starting from two skew lines and using Liaison
Addition, as discussed in \cite{M}, Section 3.3. It is easy to show (by induction on $n$) that the minimal
free resolution of $C_n$ is $$0\lra S(-2n-2)^n\lra S(-2n-1)^{4n}\lra S(-2n)^{3n+1}\lra I_{C_n}\lra 0.$$
Analogously, since the minimal free resolution of the general plane section $X_1$ of $C_1=$ two skew
lines is $$0\lra R(-3)\lra R(-2)\oplus R(-1)\lra I_{X_1}\lra 0,$$ using the short exact sequence (see
\cite{Sc}, or \cite{M}, Section 3.2 for a description of Liaison Addition and details on these
techniques) 
$$0\lra R(-2n)\lra I_{X_1}(-2n+2)\oplus I_{X_{n-1}}(-2)\lra I_{X_n}\lra 0$$ we can compute the minimal
free resolution of $I_{X_n}$, that turns out to be $$0\lra R(-2n-1)^n\lra R(-2n)\oplus R(-2n+1)^n \lra
I_{X_n}\lra 0.$$ Therefore, the degree matrix of the general plane section $X_n$ of $C_n$ is 
$$\begin{array}{ccc} 
\left(\begin{array}{c} 
1 \\ \vdots \\ 1 \end{array}\right. 
&
\underbrace{\left.\begin{array}{ccc}
2 & \cdots & 2 \\
\vdots & & \vdots \\
2 & \cdots & 2 
\end{array}\right)}_n
& \left. \begin{array}{c} \\ \\ \\ \end{array}\right\}\end{array} n$$
In this series of examples, $\M_C=soc\; H^1_*(\I_X)$, so equality is attained in Proposition~\ref{Buch},
Proposition~\ref{BuchBounds} and Corollary~\ref{dimMB}. 

Only one of the minimal generators of $I_X$ lifts to $I_C$: the one of maximum degree
$d_1$, corresponding to $a_{1,1}=1$, as shown in Corollary~\ref{liftB}.
\end{ex}

\begin{rem}
Notice that if $a_{1,t+1}=2$, that is the case for generic points in $\PP^2$ whose degree $d$ is not a
binomial coefficient ($d\neq \bin{n}{2}$ for all $n$), the deficiency module $\M_C$ has to be
concentrated in degree $a_{1,1}+\ldots +a_{t,t}-1$.

In particular, all the minimal generators of $I_X$ that are not in the initial degree, lift to $I_C$.
This also follows from the well known fact that $\alpha(I_X)\leq\alpha(I_C)\leq \alpha(I_X)+1$ (see
\cite{GM}, Corollary 3.9).
\end{rem}

So we have the following easy consequence.

\begin{cor}
Let $C\subset\PP^3$ be an arithmetically Buchsbaum curve of degree $d\neq \bin{n}{2}$ for all $n$. Assume
that the general plane section of $C$ consists of generic points.
Then the deficiency module of $C$ has to be concentrated in one degree.
\end{cor}

We now show that the bounds on the dimension of $\M_C$ of Proposition~\ref{BuchBounds} and 
Proposition~\ref{dimMB} are sharp, at least for the case of curves in $\PP^3$ and points in $\PP^2$.

\begin{thm}\label{sh}
Let $M$ be a degree matrix with at least one entry equal to $2$.
Then there exists an arithmetically Buchsbaum curve $C\subset\PP^3$ whose general plane section has degree
matrix $M$, and such that the dimension of the deficiency module $\M_C$ in each degree achieves the bound 
of Proposition~\ref{dimMB}. Moreover, $\M_C$ achieves the bounds for the initial and final degree of 
Proposition~\ref{BuchBounds}.
\end{thm}

\begin{proof}
In Remark~\ref{blocks}, we noticed that the number of non-zero components of the deficiency module is 
bounded above by the number of blocks of 2's in the degree matrix $M$. Notice that if the dimension
of $\M_C$ as a $k$-vector space is the maximum possible, according to Proposition~\ref{dimMB}, then 
the dimension of $(\M_C)_i$ for each $i$ is the maximum possible. Moreover, in this situation, all the 
graded components that
can possibly be non-zero are different from zero. Hence the bounds of Proposition~\ref{BuchBounds} on the 
initial and final degree of $\M_C$ are also attained. Therefore, in order to prove that the 
bounds of Proposition~\ref{dimMB} in every degree and the bounds of Proposition~\ref{BuchBounds} are sharp, 
it is enough to construct a curve whose deficiency module has maximum possible dimension globally.
We indicate the maximum possible dimension for $\M_C$ by $\delta(M)$, since it depends on the 
entries of the degree matrix $M$.
We prove the thesis by induction on $\delta(M)$. 
Following the notation of Proposition~\ref{dimMB}, we let
$$J=\{\: j \; | \; d_j=m_{k(j)}-n\;\; \mbox{for some $k(j)$}\}$$
and 
$$\delta(M)=\sum_{j\in J} \mmin \{ \lambda(j), \mu(j) \}.$$
Here $\lambda(j)$ is the number of rows that equal the $k(j)$-th row,
and $\mu(j)$ is the number of columns that equal the $j$-th column (the entries on the intersection
of these rows and columns form a block of 2's inside $M$, by our choice of $k(j)$). 

If $\delta(M)=1$, we can let $C$ be the curve that we constructed in Theorem~\ref{mainBuch}. 
These curves are all in the linkage class of two skew lines, hence they have $\delta(M)=1$.

So assume that we know the thesis for $\delta(M)-1$, and prove it for $\delta(M)$.
Let $(i,j)$ be such that $a_{i,j}=2$, and assume that $j\leq i$. 
Let $N$ be the submatrix of $M$ obtained by deleting the $i$-th row 
and the $j$-th column of $M$. The entries on the diagonal of $N$ are 
$a_{1,1},\ldots,a_{j-1,j-1},a_{j,j+1},\ldots,a_{i-1,i},a_{i+1,i+1},\ldots,a_{t,t}$. 
They are positive, so $N$ is a degree matrix with $\delta(N)=\delta(M)-1$. By the induction hypothesis
we have an arithmetically Buchsbaum curve $D$ with $dim(\M_D)=\delta(N)$, whose general plane section
$Y$ has degree matrix $N$.
Let $E$ be two skew lines. Let $Z=CI(1,2)$ be a general plane section of $E$.
Using Liaison Addition, we look at $I_C=FI_E+QI_D$ where $Q$ is a minimal generator of $I_E$ and 
$F$ is a form of degree 
$$a=a_{1,1}+\ldots+a_{j-1,j-1}+a_{j,j+1}+\ldots+a_{i-1,i}+a_{i+1,i+1}+\ldots+a_{t,t}-1+a_{i,t+1}$$
in the ideal of $I_D$. Notice that 
$$a-(a_{1,1}+\ldots+a_{j-1,j-1}+a_{j,j+1}+\ldots+a_{i-1,i}+a_{i+1,i+1}+\ldots+a_{t,t}+1)=$$
$$a_{i,t+1}-2\geq 0.$$
Therefore
$$\alpha(I_D)\leq a_{1,1}+\ldots+a_{j-1,j-1}+a_{j,j+1}+\ldots+a_{i-1,i}+a_{i+1,i+1}+\ldots+a_{t,t}+1\leq a,$$
and we can find a form $F$ as claimed.
By Theorem 3.2.3 in \cite{M} we have that:
\begin{itemize}
\item as sets, $C=D\cup E\cup CI(2,a)$ and
\item $\M_C\cong \M_D(-2)\oplus\M_E(-a).$
\end{itemize}
In particular, $C$ is an arithmetically Buchsbaum curve and $$dim(\M_C)=\delta(N)+1=\delta(M).$$
We still need to prove that the general plane section of $C$ has degree matrix $M$.
Let $X$ be a general plane section of $C$. Then $I_X=FI_Z+QI_Y$, and we have the short exact sequence
$$0\lra R(-a-2)\lra I_Y(-2)\oplus I_Z(-a)\lra I_X\lra 0.$$
Using the Mapping Cone argument, we obtain a free resolution for $I_X$ of the form
\begin{equation}\label{reX}
\begin{array}{rcccl}
 & R(-2-a)\oplus R(-3-a) & & \FF_1(-2) & \\
0\lra & \oplus & \lra & \oplus &\lra I_X\lra 0 \\
 & \FF_2(-2) & & R(-2-a)\oplus R(-1-a) &
\end{array}\end{equation}
where $$0\lra \FF_2\lra\FF_1\lra I_Y\lra 0$$ is a minimal free resolution for $I_Y$.
Since the image of $Q$ in $I_Z$ is a minimal generator, the free summands $R(-2-a)$ cancel in (\ref{reX}).
No other cancellation can take place, because all the other free summands come from the same minimal free 
resolution (the one of $I_Y(-2)\oplus I_Z(-a)$), so the maps between them are left unchanged under the 
Mapping Cone.
Then $X$ has minimal free resolution
$$0\lra R(-3-a)\oplus\FF_2(-2) \lra R(-1-a)\oplus\FF_1(-2)\lra I_X\lra 0$$ and its degree matrix 
has size $t\times(t+1)$, and entries as follows.
$N$ is a submatrix of it, coming from the submap $\FF_2(-2)\lra \FF_1(-2)$.
To obtain the degree matrix of $X$ from $N$, we add a row and a column corresponding to the map 
$R(-3-a)\lra R(-1-a)\oplus\FF_1(-2)$ for the row, and $R(-3-a)\oplus\FF_2(-2)\lra R(-1-a)$ for the column.
Then the entry in the intersection between the row and the column is $3+a-(1+a)=2$.
By homogeneity, the other entries on the row and column that we are adding are determined by only one of 
them. For example, the highest entry in the row is
$$3+a-(a_{1,1}+\ldots+a_{j-1,j-1}+a_{j,j+1}+\ldots+a_{i-1,i}+a_{i+1,i+1}+\ldots+a_{t,t}+2)=a_{i,t+1},$$
that coincides with the highest entry in the $i$-th row of $M$.
This proves that the degree matrix of $X$ is $M$.

We now examine the case when $a_{i,j}=2$, for some $j>i$. Pick the maximum $i$ and the minimum $j$ for 
which $a_{i,j}=2$. We can also assume that $a_{k,l}\neq 2$ for $k\leq l$.
Proceed by induction on the size $t$ of $M$. If $t=1$ the only possibility is $M=(1,2)$ and we can take
$C$ to be two skew lines. Consider the matrix $M$ of size $t\times(t+1)$, and let $N$ be the submatrix 
consisting of the last $t-1$ rows and last $t$ columns
$$N=\left(\begin{array}{cccc}
a_{2,2} & a_{2,3} & \cdots & a_{2,t+1} \\
\vdots & \vdots & & \vdots \\
a_{t,2} & a_{t,3} & \cdots & a_{t,t+1}
\end{array}\right).$$
Let $D$ be an arithmetically Buchsbaum curve, whose general plane section $Y$ has degree matrix $N$ and 
whose deficiency module has dimension $\delta(N)$. The induction hypothesis on $t$ gives the existence of 
$D$. If $\delta(N)=\delta(M)$, let $S$ be a surface of degree $a_{1,2}+\ldots+a_{t,t+1}$ containing $D$.
Such an $S$ exists since $a_{1,2}+\ldots+a_{t,t+1}\geq 1+a_{2,2}+\ldots+a_{t,t}\geq \alpha(I_D)$, 
by Proposition~\ref{gmm}. Let $T$ be a generic surface of degree $a_{1,1}$. Then $C=D\cup(S\cap T)$ is 
bilinked to $D$, therefore $dim(\M_C)=dim(\M_D)=\delta(M)$. The general plane section of $C$ has degree
matrix $M$, by Proposition 5.2.10 in \cite{M}. No cancellation occurs by genericity of the choice of $T$.
If $\delta(N)=\delta(M)-1$, then we can let $a_{i,j}=a_{1,j}=2$ for some $j\geq 2$.
By the induction hypothesis, we have an arithmetically Buchsbaum curve $D$, whose general plane section $Y$
has degree matrix $N$, and such that $dim(\M_D)=\delta(N)$. Let $E$ be a curve in the linkage class of 
two skew lines with general plane section $Z=CI(2,a_{1,2})$. Existence of $E$ follows from 
Theorem~\ref{mainBuch}.
Using Liaison Addition, let $I_C=FI_E+GI_D$ where $G$ is an element of $I_E$ of degree 2 and 
$F$ is a form of degree $a=a_{2,3}+\ldots+a_{t,t+1}$ in the ideal of $I_D$. 
$F$ can be chosen such that its image in $I_Y$ is a minimal generator, since the first column of $N$ 
has no entry equal to $2$ (see also Corollary~\ref{liftB}).
By Theorem 3.2.3 in \cite{M} we have that:
\begin{itemize}
\item as sets, $C=D\cup E\cup CI(2,a)$ and
\item $\M_C\cong \M_D(-2)\oplus\M_E(-a).$
\end{itemize}
In particular, $C$ is an arithmetically Buchsbaum curve and $$dim(\M_C)=\delta(N)+1=\delta(M).$$
We still need to prove that the general plane section of $C$ has degree matrix $M$.
Let $X$ be a general plane section of $C$. Then $I_X=FI_Z+GI_Y$, and we have the short exact sequence
$$0\lra R(-a-2)\lra I_Y(-2)\oplus I_Z(-a)\lra I_X\lra 0.$$
Using the Mapping Cone argument, we obtain a free resolution for $I_X$ of the form
\begin{equation}\label{re2X}
\begin{array}{rcccl}
 & R(-a-2)\oplus R(-2-a_{2,2}-a) & & \FF_1(-2) & \\
0\lra & \oplus & \lra & \oplus &\lra I_X\lra 0 \\
 & \FF_2(-2) & & R(-a_{2,2}-a)\oplus R(-2-a) &
\end{array}\end{equation}
where $$0\lra \FF_2\lra\FF_1\lra I_Y\lra 0$$ is a minimal free resolution for $I_Y$.
Since the image of $F$ in $I_Y$ is a minimal generator, the free summand $R(-2-a)$ cancels
with a free summand of $\FF(-2)$ in (\ref{re2X}).
No other cancellation can take place, because all the other free summands come from the same minimal free 
resolution (the one of $I_Y(-2)\oplus I_Z(-a)$), so the maps between them are left unchanged under the 
Mapping Cone. Let $\FF_1=\FF'_1\oplus R(-a)$.
Then $X$ has minimal free resolution
$$0\lra R(-2-a_{1,2}-a)\oplus\FF_2(-2) \lra R(-a_{1,2}-a)\oplus R(-2-a)\oplus\FF'_1(-2)\lra I_X\lra 0.$$
The degree matrix of $X$
has size $t\times(t+1)$, and entries as follows.
The last $t-1$ columns of $N$ are contained in it, since they come from the submap 
$\FF_2(-2)\lra \FF'_1(-2)$.
To obtain the degree matrix of $X$ from this, we add a row and two columns corresponding to the maps 
$R(-2-a_{1,2}-a)\lra R(-a_{1,2}-a)\oplus R(-2-a)\oplus\FF'_1(-2)$ for the row, and 
$R(-2-a_{1,2}-a)\oplus\FF_2(-2)\lra R(-a_{1,2}-a)\oplus R(-2-a)$ for the column.
Then the entries in intersection between the row and the columns are $2+a_{1,2}+a-a_{1,2}-a=2$, and 
$2+a_{1,2}+a-2-a=a_{1,2}$.
By homogeneity, the other entries on the row and columns that we are adding are determined by only one of 
them. For example, the highest entry in the row is
$$2+a_{1,2}+a-(a_{2,2}+\ldots+a_{t,t}+2)=a_{1,2}+a_{3,3}+\ldots+a_{t,t}+a_{2,t+1}-(a_{2,2}+\ldots+a_{t,t})=
a_{1,t+1},$$
that coincides with the highest entry in the first row of $M$.
This proves that the degree matrix of $X$ is $M$.
\end{proof}

\begin{rem}
For each 
degree matrix $M$ containing at least a $2$, one can construct an arithmetically Buchsbaum curve $C$ whose 
general plane section has degree matrix $M$ and whose deficiency module has dimension $d$ for each 
$1\leq d\leq \delta(M)$.
One way to do that is to start from the curves that we constructed in Theorem~\ref{mainBuch}, and then use 
Liaison Addition (possibly more than once) in an analogous way to what was done in the proof of 
Theorem~\ref{sh}.
\end{rem}

\end{document}